\documentclass[12 pt]{amsart}

\usepackage [mathscr]{eucal}
\usepackage{amsmath,amsthm,amssymb,amscd}
\usepackage{amsrefs}
\usepackage{epsf}
\usepackage{amsfonts}
\usepackage{dsfont}
\usepackage{latexsym}
\usepackage{mathrsfs}
\usepackage{layout}
\usepackage{bm}

\newtheorem{theorem}{Theorem}[section]
\newtheorem{lemma}[theorem]{Lemma}

\newtheorem{definition}[theorem]{Definition}
\numberwithin{equation}{section}

\DeclareMathOperator *{\essosc}{ess\ osc}
\DeclareMathOperator *{\osc}{osc}
\DeclareMathOperator *{\esssup}{ess\ sup}
\DeclareMathOperator *{\essinf}{ess\ inf}
\DeclareMathOperator *{\di}{div} 
\DeclareMathOperator *{\meas}{meas}
\DeclareMathOperator *{\dist}{dist}
\DeclareMathOperator *{\data}{data}
\DeclareMathOperator *{\diam}{diam}
\DeclareMathOperator *{\loc}{loc}
\DeclareMathOperator *{\Lip}{Lip}
\DeclareMathOperator *{\BMO}{BMO}
\DeclareMathOperator *{\BBR}{\mathbb{R}}
\DeclareMathOperator *{\BBC}{\mathbb{C}}
\newcommand{\Z}{{\mathbb Z}}

\begin{document}
\title[BMO Dirichlet problem and $A_\infty$]{BMO solvability and the $A_\infty$ condition for second order parabolic operators}

\author{Martin Dindo\v{s}}
\address{School of Mathematics, \\
         The University of Edinburgh and Maxwell Institute of Mathematical Sciences}
\email{M.Dindos@ed.ac.uk}

\author{Stefanie Petermichl}
\address{Institut de Math\'ematiques de Toulouse}
\email{stefanie.petermichl@math.univ-toulouse.fr}

\author{Jill Pipher}
\address{Department of Mathematics, Brown University}
\email{Jill.Pipher@Brown.edu}


\begin{abstract} We prove that the $A_\infty$ property of parabolic measure for operators in certain time-varying domains is equivalent to a Carleson measure property of bounded solutions.
In \cite{KKKT}, this criterion on bounded solutions was established in the elliptic case, improving an earlier result in \cite{DKP} for solutions with data in $BMO$.
The extension to the parabolic setting requires
an approach to the key estimate of \cite{KKPT} and
\cite{KKKT} that primarily exploits the maximum principle. For various classes of parabolic operators (\cite{Rn2}), this criterion makes it easier to establish the solvability of the Dirichlet problem with data in $L^p$ for some $p$, and also to quantify these results in several aspects. 
\end{abstract}

\maketitle
\section{Introduction}\label{S0:Intro}

In this paper, we prove a criterion for establishing the $L^p$ solvability of the Dirichlet problem for parabolic operators $L=\partial_t-\di (A \nabla \cdot)$ in certain time-varying domains, and where the
matrix $A$ satisfies an ellipticity condition. 
Our results are analogous to similar criteria established in \cite{DKP} and \cite{KKKT} in for elliptic operators  $\di (A \nabla \cdot)$.
With this criterion we are then able to give a simpler proof of existence of $L^p$ solvability, for some $p$, for a class of operators, studied in \cite{Rn}, whose coefficients satisfy a Carleson-measure regularity condition, also permitting us to quantify the dependence on $p$.
(See Section 6.)
In virtue of the maximum principle, there is a natural (``parabolic") representing measure associated with the solvability of the Dirichlet problem for such $L$ with continuous data. In the domains considered here, this measure has been shown to be doubling in \cite{HL}. 
We will be interested in $L^p$ solvability of 
boundary value problems with respect to a natural measure $\sigma$ (see Definition \ref{D:measure}) defined on the boundary of the time-varying domain, one which coincides with surface measure when that domain
is sufficiently smooth.

The study of the heat equation in non-smooth domains, or more generally of parabolic operators with non-smooth coefficients, has historically closely followed the development of the elliptic theory, while presenting new challenges to finding the correct analogues of the elliptic results.

Dahlberg \cite{Da} showed that, in a Lipschitz domain, harmonic measure and surface measure are mutually absolutely continuous,
and that in fact the elliptic Dirichlet problem is solvable with data in $L^{2}$ with respect to surface measure. R. Hunt then asked whether
Dalhberg's result held for the heat equation in
domains whose boundaries are given locally as functions $\psi(x,t)$, Lipschitz in the spatial variable.
It was natural to conjecture that the correct regularity of $\psi(x,t) $ in the time variable $t$ should be a H\"older condition of order $1/2$ in $t$ (denoted $\Lip_{1/2}$ in $t$). However, the counterexamples of \cite{KW} showed that
this condition did not suffice. Lewis and Murray \cite{LM} then established mutual absolute
continuity of caloric measure and a certain parabolic analogue of surface measure when $\psi$ has $1/2$ of a time derivative in the parabolic $\BMO (\mathbb{R}^n)$ space, a slightly stronger condition than $\Lip_{1/2}$. Hofmann and Lewis (\cite{HL1}) subsequently solved the $L^2$ Dirichlet problem for the heat equation in 
graph domains of Lewis-Murray type when the $BMO$ norm of the time derivative was sufficiently small.


In this paper, we consider parabolic operators of the form
\begin{equation}\label{E:v}\begin{cases}
u_t = \di (A \nabla u) & \text{in } \Omega, \\
u   = f               & \text{on } \partial \Omega
\end{cases}\end{equation}
where $A= [a_{ij}(X, t)]$ is an $n\times n$ matrix satisfying a uniform ellipticity condition: there exists positive constants $\lambda$ and $\Lambda$ such that
\begin{equation}\label{E:UEBv}
\lambda |\xi|^2 < \sum_{i,j} a_{ij} \xi_i \xi_j < \Lambda |\xi|^2
\end{equation}
for all $\xi \in \mathbb{R}^{n}$. 

Here and throughout we will consistently use $\nabla u$ to denote the gradient in the spatial variables, $u_t$ or $\partial_t u$ the gradient in the time variable and use $Du=(\nabla u, \partial_t u)$ for the full gradient of $u$.\vglue2mm
As in  \cite{DH}, the results here are formulated for the class of {\it admissible parabolic domains}, which are, in effect,  bounded time-varying domains that are ``locally" of Lewis-Murray type. A related, but smaller, class of localized domains in which parabolic boundary value problems are solvable was considered in \cite{Rn2}.

It is a fact that the parabolic PDE (\ref{E:para}) with continuous boundary data is uniquely solvable (c.f. discussion under Definition 2.7 in \cite{DH}) and that there exists a a measure $\omega^{(X,t)}$ such that
\begin{equation}u(X,t)=\int_{\partial \Omega} f(y,s) d\omega^{(X,t)}(y,s)\label{E:pm}\end{equation}
for all continuous data, called the parabolic measure. Under the assumption of Definition \ref{D:domain}, this measure is doubling (c.f. \cite{N}).

As $\omega^{(X,t)}$ is a Borel measure, if follows that we can use \eqref{E:pm} to extend the solvability to a a class of bounded Borel measurable functions $f$. This observation will be important later.
\vglue1mm

In Sections 2 and 3, we recall the definitions of parabolic measure and the fact that solvability of the Dirichlet problem for $L$ with data in some $L^p$ space, for some $p<\infty$ is equivalent to 
the $A_\infty$ property relative to the boundary measure $\sigma$.  Our main results are the following.

\begin{theorem}\label{T:Main} Let $\Omega$ be a domain as in Definition \ref{D:domain} with character $(\ell,N,C_0)$. Let $A=[a_{ij}]$ be a matrix with bounded
measurable coefficients defined on $\Omega$ satisfying uniform ellipticity and boundedness with constants $\lambda$ and $\Lambda$.

\noindent If the parabolic measure for the operator $L=\partial_t-\di (A \nabla \cdot)$ is in the class $A_{\infty}(d\sigma)$, then the BMO Dirichlet boundary
value problem defined in Definition \ref{D:bmo} is solvable and the estimate
\[\sup_{\Delta\subset\partial\Omega}\sigma(\Delta)^{-1}\int_{T(\Delta)}|\nabla u|^2\delta\,dXdt \lesssim \|f\|^2_{BMO(\partial \Omega, d\sigma)}
\]
holds for all continuous functions $f\in C(\partial\Omega)$. Here $\delta = \delta(X,t)$ is the parabolic distance to the boundary.
\end{theorem}

\begin{theorem}\label{T:Main2} Let $\Omega$ be a domain as in Definition \ref{D:domain} with character $(\ell,N,C_0)$. Let $A=[a_{ij}]$ be a matrix with bounded
measurable coefficients defined on $\Omega$ satisfying uniform ellipticity and boundedness with constants $\lambda$ and $\Lambda$.

\noindent Assume that for all continuous functions $f\in C(\partial\Omega)$ the corresponding solution $u$ satisfies the estimate
\begin{equation}
\sup_{\Delta\subset\partial\Omega}\sigma(\Delta)^{-1}\int_{T(\Delta)}|\nabla u|^2\delta\,dX \le C \|f\|^2_{L^\infty(\partial \Omega, d\sigma)},\label{E:bmoinfty}
\end{equation}
for a constant $C=C(\Omega,A)>0$. Then the parabolic measure for the operator $L=\partial_t-\di (A \nabla \cdot)$ belongs to the class $A_{\infty}(d\sigma)$.
Hence for some $p_0<\infty$ the $L^p$ Dirichlet boundary value problem for the operator $L$ is solvable for all $p\in (p_0,\infty)$.
\end{theorem}

Both of these theorems are parabolic analogues of established results for elliptic operators: see \cite{DKP} and \cite{KKKT}.  Our proof of Theorem \eqref{T:Main2} uses the 
primary strategy laid out in \cite{KKPT} and \cite{KKKT}, but with a simpler approach to the key estimate in order to adapt it to the parabolic setting.

A key feature of Theorem \ref{T:Main2} is that one only needs to check the bound \eqref{E:bmoinfty} using the $L^\infty$ norm of $f$, as opposed to a BMO norm. This condition is easier to verify since the BMO norm of a function can be smaller than its $L^\infty$ norm. The analogous elliptic result with $L^\infty$ norm,  stated below, was established in \cite{KKKT} 
and the proof presented here also easily goes over to the elliptic setting.

\begin{theorem}\label{T:Main3} (\cite{KKKT}) Let $\Omega$ be a Lipschitz domain and $A$ be a uniformly elliptic matrix on $\Omega$ with bounded measurable coefficients.
If for all continuous functions $f\in C(\partial\Omega)$ the corresponding solution $u$ of the equation $\di (A \nabla u)=0$ satisfies the estimate
\begin{equation}
\sup_{\Delta\subset\partial\Omega}\sigma(\Delta)^{-1}\int_{T(\Delta)}|\nabla u|^2\delta\,dX \le C \|f\|^2_{L^\infty(\partial \Omega, d\sigma)},\label{E:bmoinfty2}
\end{equation}
for a constant $C=C(\Omega,A)>0$, then the elliptic measure for the operator $L=\di (A \nabla \cdot)$ belongs to the class $A_{\infty}(d\sigma)$.
Hence for some $p_0<\infty$ the elliptic $L^p$ Dirichlet boundary value problem for the operator $L$ is solvable for all $p\in (p_0,\infty)$.
\end{theorem}

Our Theorem \ref{T:Main2} also provides an easier proof of the main result of \cite{Rn2} on parabolic operators with coefficients whose gradients satisfy Carleson condition, or a slightly weaker assumption on the oscillation. This complements the results of 
 \cite{DH} where it was established that if 
\begin{equation}\label{carlI}
d\mu=\delta(X)^{-1}\left(\mbox{osc}_{B_{\delta(X)/2}(X)}a_{ij}\right)^2\,dX\,dt
\end{equation}
is a density of small Carleson measure with norm $\|\mu\|_{Carl}$ on an admissible parabolic domain $\Omega$, then given $2\le p<\infty$ there exists a constant $C(p,\lambda,\Lambda)>0$ such that for $\max\{\ell^2,\|\mu\|_{Carl}\}<C(p)$ then the $L^p$ Dirichlet problem for the operator $L=\partial_t-\di (A \nabla \cdot)$ is solvable.

The next result, with no smallness assumptions, is a quantitative version of Theorem 3.3 of \cite{Rn2}. 

\begin{theorem}\label{T:Main4} (A Quantitative version of Theorem 3.3 of \cite{Rn2}.)  
 Let $\Omega$ with character $(\ell,N,C_0)$ and $A$ be as in Theorem \ref{T:Main}. Denote by $\mu$ the measure with density
\begin{equation}\label{carlII}
d\mu=\delta(X)^{-1}\left(\mbox{osc}_{B_{\delta(X)/2}(X)}a_{ij}\right)^2\,dX\,dt
\end{equation}
and by $\|\mu\|_{Carl}$ its Carleson norm.

\noindent For every $2\le p\le \infty$ there exists a constant $C(p,\lambda,\Lambda)>0$ such that for $\max\{\ell^2,\|\mu\|_{Carl}\}<C(p)$ then the $L^p$ Dirichlet problem for the operator $L=\partial_t-\di (A \nabla \cdot)$ is solvable. Moreover, for $\lambda,\Lambda$ fixed the constant
$$C(p,\lambda,\Lambda)\to \infty,\qquad\mbox{as}\quad p\to\infty.$$

It follows that for any admissible domain $\Omega$ with character $(\ell,N,C_0)$ and any parabolic operator $L=\partial_t-\di (A \nabla \cdot)$ with ellipticity
constants $\lambda, \Lambda$, if the $\mu$ defined by \eqref{carlII} satisfies $\ell<\infty,\,\|\mu\|_{Carl}<\infty$, then there exists
$$p_0=p_0(\lambda,\Lambda,\ell,n,\|\mu\|_{Carl}) <\infty$$
such that for all $p\in (p_0,\infty)$  the $L^p$ Dirichlet problem for the operator $L$ is solvable. In particular, the parabolic measure of $L$ belongs to $B_{p'}(d\sigma)\subset A_{\infty}(d\sigma)$ ($p'=p/(p-1)$).
\end{theorem}

In \cite{Rn2}, it has been shown that there exists a $p$ for which the $L^p$ Dirichlet problem is solvable when $\|\mu\|_{Carl}$  is finite.
However, it would not, given the method of proof, be easy to track the dependence of $p$ on the Carleson norm, nor would it be possible to address solvability for a particular value of $p$. 
The fact that for every $2\le p\le \infty$ there exists a constant $C(p,\lambda,\Lambda)>0$ such that the $L^p$ Dirichlet problem for the operator $L=\partial_t-\di (A \nabla \cdot)$ is solvable whenever
$\max\{\ell^2,\|\mu\|_{Carl}\}<C(p)$ is in \cite{DH}. But it was not evident that $C(p)=C(p,\lambda,\Lambda)\to \infty \,\mbox{as}\, p\to\infty$. In fact, the estimates from below for $C(p)$ in \cite{DH} are not (due to the method employed) powerful enough to show that $C(p)\to\infty$ as $p\to\infty$. This, then, is the main contribution of Theorem \ref{T:Main4}.

\section{Preliminaries}\label{S1:Pre}
\subsection{Admissible parabolic domains.}

In this subsection we recall the class of ``admissible" time-varying domains in \cite{DH} whose boundaries are given locally as functions $\psi(x,t)$, Lipschitz in the spatial variable and
satisfying the Lewis-Murray condition in the time variable. At each time $\tau\in\BBR$ the set of points in $\Omega$ with fixed time $t=\tau$, that is
$\Omega_\tau=\{(X,\tau)\in\Omega\}$ will be assumed to be a nonempty bounded Lipschitz domain in $\BBR^n$.  We choose to consider domains that are bounded (in space) since
this most closely corresponds to domains considered in the paper \cite{DPP} (for the elliptic equation).
However, our result can be adapted to the case of unbounded domains (in space). See \cite{HL} which focuses on the unbounded case.

We start with few preliminary definitions, formulated exactly as in \cite{DH}.

If
 $\psi(x,t): \mathbb{R}^{n-1} \times \mathbb{R} \to \mathbb{R}$ is a compactly supported function,  the half derivative in time may be defined by the Fourier transform or by
\[
D_{1/2}^{t} \psi(x,t) = c_n \int_{\mathbb{R}} \frac{\psi(x,s) - \psi(x,t)}{|s-t|^{3/2}} \,ds
\]
for a properly chosen constant $c_n$ (depending on the dimension $n$).

We shall also need a local version of this definition. If $I\subset \BBR$ is a bounded interval and $\psi(x,t)$ is defined on $\{x\}\times I$ we consider:
\[
D_{1/2}^{t} \psi(x,t) = c_n \int_{I} \frac{\psi(x,s) - \psi(x,t)}{|s-t|^{3/2}} \,ds,\qquad\mbox{for all }t\in I.
\]

We define a parabolic cube in $\mathbb{R}^{n-1} \times \mathbb{R}$,  for a constant $r>0$, as
\begin{equation}\label{D:Q}
Q_{r} (x, t)
= \{ (y, s) \in \mathbb{R}^{n-1}\times\mathbb{R} : |x_i - y_i| < r \ \text{for all } 1 \leq i \leq n-1, \ | t - s |^{1/2} < r \}.
\end{equation}
For a given $f: \mathbb{R}^{n} \to \mathbb{R}$ let,
\[
f_{Q_r} = |Q_r|^{-1} \int_{Q_r} f(x,t) \,dx\,dt.
\]
When we write $f \in \BMO (\mathbb{R}^n) $ we mean that $f$ belongs to the parabolic version of the usual BMO space with the norm $\|f\|_{*}$ where
\[
\|f\|_{*} = \sup_{Q_r} \left\{  \frac{1}{|Q_r|} \int_{Q_r} |f - f_{Q_r} | \,dx\,dt \right\} < \infty.
\]

Again, we also consider a local version of this definition. For a function $f:J\times I\to\BBR$, where $J\subset \BBR^{n-1}$ and $I\subset \BBR$ are closed bounded balls we consider the norm $\|f\|_{*}$ defined as above where the supremum is taken over all parabolic cubes $Q_r$ contained in $J\times I$.\vglue3mm

The following definitions are motivated by the standard definition of a Lipschitz domain.

\begin{definition}
$\Z \subset \BBR^n\times \BBR$ is an $\ell$-cylinder of diameter $d$ if there
exists a coordinate system $(x_0,x,t)\in \BBR\times\BBR^{n-1}\times \BBR$ obtained from the original coordinate system only by translation in spatial and time variables and rotation in the spatial variables such that
\[
\Z = \{ (x_0,x,t)\; : \; |x|\leq d,\; |t|\leq d^2, \; -(\ell+1)d \leq x_0 \leq (\ell+1)d \}
\]
and for $s>0$,
\[
s\Z:=\{(x_0,x,t)\;:\; |x|<sd,\; |t|\leq s^2d^2, \; -(\ell+1)sd \leq x_0 \leq (\ell+1)sd \}.
\]
\end{definition}

\begin{definition}\label{D:domain}
$\Omega\subset \BBR^n\times \BBR$ is an {\it admissible parabolic} domain
with `character' $(\ell,N,C_0)$ if there exists a positive scale $r_0$ such that for any time $\tau\in\BBR$
there are at most $N$ $\ell$-cylinders $\{{\Z}_j\}_{j=1}^N$ of diameter $d$, with
$\frac{r_0}{C_0}\leq d \leq C_0 r_0$ such that\vglue2mm

\noindent (i) $8{\Z}_j \cap \partial\Omega$ is the graph $\{x_0=\phi_j(x,t)\}$ of a
function $\phi_j$, such that
\begin{equation}
|\phi_j(x,t)- \phi_j(y,s)| \leq \ell[|x-y|+|t-s|^{1/2}], \qquad
\phi_j(0,0)=0\label{E:L1}
\end{equation}
and
\begin{equation}
\|D^t_{1/2}\phi_j\|_*\le \ell.\label{E:L2}
\end{equation}
\vglue2mm

\noindent (ii) $\displaystyle \partial\Omega\cap\{|t-\tau|\le d^2\}=\bigcup_j ({\Z}_j \cap \partial\Omega
)$,

\noindent (iii) In the coordinate system $(x_0,x,t)$ of the $\ell$-cylinder $\Z_j$:
$$\displaystyle{\Z}_j \cap \Omega \supset \left\{
(x_0,x,t)\in\Omega \; : \; |x|<d, \;|t|<d^2\;, \delta(x_0,x,t)=\mathrm{dist}\left( (x_0,x,t),\partial\Omega
\right) \leq \frac{d}{2}\right\}.$$
\end{definition}

Here the distance is the parabolic distance $d[(X,t),(Y,\tau)]=(|X-Y|^2+|t-\tau|)^{1/2}$.\vglue2mm

\noindent{\it Remark.} It follows from this definition that for each $\tau\in\BBR$ the time-slice $\Omega_\tau=\Omega\cap\{t=\tau\}$ of an admissible parabolic domain $\Omega\subset \BBR^n\times \BBR$ is a bounded Lipschitz domain in $\BBR^n$ with `character' $(\ell,N,C_0)$. Due to this fact, the Lipschitz domains $\Omega_{\tau}$ for all $\tau\in\BBR$ have all uniformly bounded diameter (from below and above). That is
$$\inf_{\tau\in\mathbb R}\diam(\Omega_\tau)\approx r_0 \approx \sup_{\tau\in\mathbb R}\diam(\Omega_\tau),$$
where $r_0$ is scale from Definition \ref{D:domain} and the implied constants in the estimate above only depend on $N$ and $C_0$.

In particular, if ${\mathcal O}\subset \BBR^n$ is a bounded Lipschitz domain, then the parabolic cylinder $\Omega={\mathcal O}\times \BBR$ is an example of a domain satisfying Definition \ref{D:domain}.\vglue2mm


\begin{definition}\label{D:measure}
Let $\Omega\subset \BBR^n\times \BBR$ be an {\it admissible parabolic} domain
with `character' $(\ell,N,C_0)$. The measure $\sigma$, defined on sets $A\subset \partial\Omega$ is:
\begin{equation}\label{E:sigma}
\sigma(A)=\int_{-\infty}^\infty {\mathcal H}^{n-1}\left(A\cap\{(X,t)\in\partial\Omega\}\right)dt,
\end{equation}
where ${\mathcal H}^{n-1}$ is the $n-1$ dimensional Hausdorff measure on the Lipschitz boundary $\partial\Omega_t=\{(X,t)\in\partial\Omega\}$.
\end{definition}

We are going to consider solvability of the $L^p$ and BMO Dirichlet boundary value problems with respect to the measure $\sigma$. The measure $\sigma$
may not be comparable to the usual surface measure on $\partial\Omega$: in the $t$-direction the functions $\phi_j$ from the Definition \ref{D:domain}
are only half-Lipschitz and hence the standard surface measure might not be locally finite.

Our definition assures that for any $A\subset \Z_j$, where  $\Z_j$ is an $L$-cylinder we have
\begin{equation}\label{E:comp}
{\mathcal H}^{n}(A) \approx \sigma\left(\{(\phi_j(x,t),x,t):\,(x,t)\in A\}\right),
\end{equation}
where the actual constants in \eqref{E:comp} by which these measures are comparable only depend on the $\ell$ of the `character' $(\ell,N,C_0)$ of domain $\Omega$.\vglue2mm

If $\Omega$ has smoother boundary, such as Lipschitz (in all variables) or better, then the measure $\sigma$ is comparable to the usual $n$-dimensional Hausdorff measure ${\mathcal H}^{n}$. In particular, this holds for a parabolic cylinder $\Omega={\mathcal O}\times \BBR$.

\begin{definition} Let $\Omega$ be an admissible parabolic domain from Definition \ref{D:domain}.
For $(Y,s)\in\partial\Omega$, $(X,t)\in \Omega$, $r>0$, and $d$ the parabolic distance we write:
\begin{align*}
B_r(Y,s) &=\{(X,t)\in{\BBR}^{n}\times \BBR:\, d[(X,t),(Y,s)]<r\}\\
\Delta_r(Y,s) &= \partial \Omega\cap B_r(Y,s),\,\,\,\qquad T(\Delta_r) = \Omega\cap B_r(Y,s).
\end{align*}
\end{definition}

\begin{definition}\label{D:carl}
Let $T(\Delta_r)$ be the Carleson region associated to a surface
ball $\Delta_r$ in $\partial\Omega$, as defined above. A measure $\mu:\Omega\to\BBR^+$  is said to be
Carleson if there exists a constant $C=C(r_0)$ such that for all
$r\le r_0$ and all surface balls $\Delta_r$
\[\mu(T(\Delta_r))\le C \sigma (\Delta_r).\]
The best possible constant $C$ will be called the Carleson norm and shall be denoted by $\|\mu\|_{C,r_0}$. We write $\mu \in \mathcal C$.  If
$\displaystyle\lim_{r_0\to 0} \|\mu\|_{C,r_0}=0$,  we say that the
measure $\mu$ satisfies the vanishing Carleson condition and write $\mu \in \mathcal C_V$.
\end{definition}

When $\partial\Omega$ is locally given as a graph of a function $x_0=\psi(x,t)$ in the coordinate system $(x_0,x,t)$
and $\mu$ is a measure supported on $\{x_0>\psi(x,t)\}$ we can reformulate the Carleson condition locally using the parabolic cubes $Q_r$
and corresponding Carleson regions $T(Q_r)$ where
\begin{align*}
Q_{r} (y, s)
&= \{ (x, t) \in \mathbb{R}^{n-1}\times\mathbb{R} : |x_i - y_i| < r \ \text{for all } 1 \leq i \leq n-1, \ | t - s |^{1/2} < r \}\\
T(Q_r)&=\{(x_0,x, t) \in \mathbb{R}\times\mathbb{R}^{n-1}\times\mathbb{R} : \psi(x,t)<x_0<\psi(x,t)+r,\, (x,t)\in Q_{r} (y, s)\}.
\end{align*}

The Carleson condition becomes
\[\mu(T(Q_r))\le C|Q_r|=Cr^{n+1}.\]

We remark, that the corresponding Carleson norm will not be equal to the one from Definition \ref{D:carl} but these norms will be comparable. Hence the notion of
vanishing Carleson norm does not change if we take this as the definition of the Carleson norm instead of Definition \ref{D:carl}.\vglue2mm

Observe also, that the function $\delta(X,t):=\inf_{(Y,\tau)\in \partial\Omega}d[(X,t),(Y,\tau)]$ that is measuring the distance of a point $(X,t)=(x_0,x,t)\in\Omega$ to the boundary $\partial\Omega$ is comparable to $x_0-\psi(x,t)$ which in turn is comparable to $[\rho^{-1}(X,t)]_{x_0}$ (the first component of the inverse map $\rho^{-1}$).\vglue2mm

\begin{definition} (Corkscrew points) Let $\Omega$ be an admissible parabolic domain from Definition \ref{D:domain} and $r_0>0$ the scale defined there.
For any boundary ball $\Delta_r=\Delta_r(Y,s)\subset\partial\Omega$ with $0<r\lesssim r_0$ we say that a point $(X,t)\in\Omega$ is a {\bf corkscrew point}
 of the ball $\Delta_r$ if
$$t=s+2r^2,\qquad \delta(X,t)\approx r\approx d[(X,t),(Y,s)].$$
That is the point $(X,t)$ is an interior point of $\Omega$ of distance to the ball $\Delta_r$ and the boundary $\partial\Omega$ of order $r$. The point $(X,t)$
lies at the time of order $r^2$ further then the times for the ball $\Delta_r$. Finally, the implied constants in the definition above only depend on the domain
$\Omega$ but not on $r$ and the point $(Y,s)$.

Each ball of radius $0<r\lesssim r_0$ has infinitely many corkscrew points; for each ball we choose one and denote it by $V(\Delta_r)$ or if there is no confusion to which ball the corkscrew point belongs just $V_r$.
\end{definition}

\noindent{\it Remark.} Given the fact that the time slices $\Omega_\tau$ of the domain $\Omega$ are of approximately diameter $r_0$ the corkscrew points do not exists for balls of sizes $r>>r_0$.

\subsection{Parabolic Non-tangential cones and related functions} \label{S12:Not}

We proceed with the definition of parabolic non-tangential cones. We define the cones in a (local) coordinate system where $\Omega=\{(x_0,x,t):\,x_0>\psi(x,t)\}$. In particular this also applies to the upper half-space $U=\{(x_0,x,t), x_0>0\}$. We note here, that a different choice of coordinates (naturally) leads to different sets of cones, but as we shall establish the particular choice of non-tangential cones is not important as it only changes constants in the estimates for the area, square and non-tangential maximal functions defined using these cones. However the norms defined using different sets of non-tangential cones are comparable.

For a constant $a>0$, we define the parabolic non-tangential cone at a point $(x_0,x,t)\in\partial\Omega$ as follows
\begin{equation}\label{D:Gammag}
\Gamma_{a}(x_0, x, t) = \left\{(y_0, y, s)\in \Omega : |y - x| + |s-t|^{1/2} < a(y_0 - x_0), \  y_0 > x_0 \right\}.
\end{equation}
We occasionally truncate the cone $\Gamma$ at the height $r$
\begin{equation}\begin{split}\label{D:Gammagr}
&\Gamma_{a}^{r}(x_0, x, t) =\\ &\left\{(y_0, y, s)\in \Omega : |y - x| + |s-t|^{1/2} < a(y_0 - x_0), \  x_0 < y_0 < x_0 + r  \right\}.
\end{split}
\end{equation}

When working on the upper half space (domain $U$), $(0, x, t)$ is the boundary point of $\partial U$. In this case we shorten the notation and write
\begin{equation}\label{D:Gamma}
\Gamma_{a}(x,t) \qquad\mbox{instead of }\qquad \Gamma_{a}(0, x, t)
\end{equation}
and
\begin{equation}\label{D:Gammar}
\Gamma_{a}^{r}(x,t) \qquad\mbox{instead of }\qquad \Gamma^r_{a}(0, x, t).
\end{equation}

Observe that the slice of the cone $\Gamma_{a}(x_0, x, t)$ at a fixed height $h$ is the set
$$\{(y,s):\,(x_0+h,y,s)\in\Gamma_a(x_0,x,t)\}$$
which contains and is contained in a parabolic box $Q_{s}(x,t)$ of radius $s$ comparable to $h$, that is
for some constants $c_1, c_2$ depending only on the dimension $n$ and $a$ we have
\[
Q_{c_1 h}(x,t) \subset \{(y,s):\,(x_0+h,y,s)\in\Gamma_a(x_0,x,t)\} \subset Q_{c_2 h}(x,t).
\]

For a function $u : \Omega \rightarrow \mathbb{R}$, {\it the nontangential maximal function} $\partial\Omega\to \BBR$ and its truncated version at a height $r$ are defined as
\begin{equation}\label{D:NTan}\begin{split}
N_{a}(u)(x_0,x,t) &= \sup_{(y_0,y,s)\in \Gamma_{a}(x_0,x,t)} \left|u(y_0 , y, s)\right|, \\
N_{a}^{r}(u)(x_0,x,t) &= \sup_{(y_0,y,s)\in \Gamma_{a}^{r}(x_0,x,t)} \left|u(y_0 , y, s)\right|\quad\mbox{for }(x_0,x,t)\in\partial\Omega.
\end{split}\end{equation}

Now we define {\it the square function} $\partial\Omega\to \BBR$ (and its truncated version) assuming $u$ has a locally integrable distributional gradient by
\begin{equation}\begin{split}\label{D:Square}
S_{a}(u)(x_0,x,t) &= \left(\int_{\Gamma_{a}(x_0,x,t)} (y_{0}-x_0)^{-n} |\nabla u|^{2}(y_0, y, s) \,dy_0\,dy\,ds\right)^{1/2}, \\
S_{a}^{r}(u)(x_0,x,t) &= \left(\int_{\Gamma_{a}^{r}(x_0,x,t)} (y_{0}-x_0)^{-n} |\nabla u|^{2}(y_0, y, s) \,dy_0\,dy\,ds\right)^{1/2}.
\end{split}\end{equation}
Observe that on the domain $U=\{(x_0,x,t):\,x_0>0\}$
\[
\|S_a (u)\|^{2}_{L^2(\partial U)} = \int_{U} y_{0} |\nabla u|^{2}(y_0, y, s) \,dy_0\,dy\,ds.
\]\vglue2mm


\subsection{$L^p$ and BMO Solvability of the Dirichlet boundary value problem}

We are now in a position to define $L^p$ solvability of a Dirichlet problem for a parabolic operator.

\begin{definition}
\label{WeakSol} (\cite{Aron}) We say that $u$ is a weak solution to $L$  in $\Omega$ if $u, \nabla u \in L^2_{{loc}}(\Omega)$ and
$ \sup_t \| u(\cdot,t) \|_{L^2_{{loc}}(\Omega_t)} < \infty$,
and 
$$ \int_{\Omega} ( -u \phi_t + A \nabla u.\nabla \phi ) dXdt = 0$$
for all $\phi \in C^{\infty}_0(\Omega)$.
\end{definition}

\begin{definition}\label{D:Lp} Let $1<p\le\infty$  and $\Omega$ be an admissible parabolic domain from the Definition \ref{D:domain}.
Consider the parabolic Dirichlet boundary value problem

\begin{equation}\label{E:para}\begin{cases}
u_t = \di (A \nabla u)  & \text{in } \Omega, \\
u   = f \in L^p              & \text{on } \partial \Omega,\\
N(u) \in L^p(\partial\Omega,d\sigma).
\end{cases}\end{equation}
where the matrix $A= [a_{ij}(X, t)]$ satisfies the uniform ellipticity condition and $\sigma$ is the measure supported on $\partial\Omega$ defined by \eqref{E:sigma}.

We say that Dirichlet problem with data in $L^p(\partial \Omega, d\sigma)$
is solvable if the (unique) solution $u$ with continuous boundary data $f$ satisfies the estimate
\begin{equation}\label{E:estim}
\|N(u)\|_{L^p(\partial \Omega, d\sigma)} \lesssim \|f\|_{L^p(\partial \Omega, d\sigma)}.
\end{equation}
The implied constant depends only the operator $L$,
$p$, and the the triple $(L,N,C_0)$ of Definition \ref{D:domain}.
\end{definition}

The $L^p$ solvability of the Dirichlet boundary value problem for some $p<\infty$ is equivalent to the parabolic measure $\omega$ belonging to a ``parabolic $A_\infty$" class
with respect to the measure $\sigma$ on the surface $\partial\Omega$ (Theorem 6.2 in \cite{N}).
We now recall the definition of parabolic $A_\infty$.

\begin{definition}($A_\infty$ and $B_p$) \label{D:infty} Let $\Omega$ be an admissible parabolic domain from Definition \ref{D:domain}. For a ball $\Delta_d$ with radius  $d\lesssim \sup_{\tau}\diam(\Omega_\tau)$ we denote its corkscrew point by $V_d$.

We say that parabolic measure of an operator $L=\partial_t - \di (A \nabla \cdot)$ is $A_\infty(\Delta_d)$ if for every $\varepsilon>0$ there exists $\delta>0 = \delta(\varepsilon)$
such that for any ball $\Delta \subset \Delta_d$ and subset $E\subset\Delta$ we have:
$$\frac{\omega^{V_d}(E)}{\omega^{V_d}(\Delta)}<\delta \Longrightarrow \frac{\sigma(E)}{\sigma(\Delta)}<\varepsilon.$$
The measure is $A_\infty$ if it belongs to $A_\infty(\Delta_d)$ for all $\Delta_d$.
If $A_\infty$ holds then the measures $\omega^{V_d}$ and $\sigma$ are mutually absolutely continuous and hence one can write $d\omega^{V_d}=K^{V_d}d\sigma$.

For $p\in (1,\infty)$ we say that $\omega$ belongs to
belongs to the reverse-H\"older class $B_p(d\sigma)$ if $K^{V_d}$ satisfies the reverse H\"older inequality
$$\left(\sigma(\Delta)^{-1}\int_{\Delta} \left(K^{V_d}\right)^p\,d\sigma\right)^{1/p}\lesssim \sigma(\Delta)^{-1}\int_{\Delta} K^{V_d}\,d\sigma,$$
for all balls $\Delta\subset \Delta_d$.
\end{definition}

\noindent {\it Remark 1.} It can be shown that $A_\infty(d\sigma)=\bigcup_{p>1}B_p(d\sigma)$.\vglue2mm

\noindent {\it Remark 2.} It has been shown in \cite{Rn} that if the parabolic measure is $A_\infty$ with respect to the surface measure $\sigma$ then
the non-tangential maximal function and the square function of a solution are equivalent, that is for all $1<p<\infty$

\[
\int_{\partial \Omega} N^p (u) \,dx\,dt \approx \int_{\partial \Omega} S^p (u) \,dx\,dt + \int_{\partial \Omega} u^{p} \,dx\,dt
\]
See also Theorem 6.2 of \cite{DH}. Here the implied constants do not depend on the solution $u$, only on $p$, the domain and the parabolic operator. Hence if follows that if \eqref{E:estim} holds then also
\[\|S(u)\|_{L^p(\partial \Omega, d\sigma)} \lesssim \|f\|_{L^p(\partial \Omega, d\sigma)}.
\]
It turns out that is condition is more convenient to define the end-point BMO Dirichlet boundary value problem.

\begin{definition}\label{D:bmo} Let $\Omega$ and the matrix $A$ be as in Definition \ref{D:Lp}.
We say that the Dirichlet problem with data in $BMO(\partial \Omega, d\sigma)$
is solvable if the (unique) solution $u$ with continuous boundary data $f$ satisfies the estimate
\begin{equation}\label{E:bmo}
\sup_{\Delta\subset\partial\Omega}\sigma(\Delta)^{-1}\int_{T(\Delta)}|\nabla u|^2\delta\,dX \lesssim \|f\|^2_{BMO(\partial \Omega, d\sigma)}.
\end{equation}
The implied constant depends only the operator $L$ and the the triple $(L,N,C_0)$ of Definition \ref{D:domain}. Here the supremum on the right-hand side
is taken over all parabolic balls $\Delta\subset\partial\Omega$. $T(\Delta)$ denotes the corresponding Carleson region (as defined above).
\end{definition}

\noindent {\it Remark 3.} The term on left-hand side of \eqref{E:bmo} is connected with the square function in the following way. If $\Delta=\Delta_r$
is a parabolic boundary ball, then
\[\sigma(\Delta)^{-1}\int_{T(\Delta)}|\nabla u|^2\delta\,dX \approx \sigma(\Delta)^{-1}\int_{\Delta}(S_r(u))^2\,d\sigma,\]
where $S_r$ is the truncated square function at height $r$. To be completely correct, in the inequalities implied by the previous line in the bounds from above we should enlarge $\Delta_r$ to its double, say $\Delta_{2r}$, this however makes no difference if we want to replace the left-hand side of \eqref{E:bmo} by the
integral over the square function as we are taking the supremum over all boundary balls $\Delta$ anyway.
\vglue1mm
\noindent {\it Remark 4.} Is is sufficient to assume \eqref{E:bmo} only holds for all balls $\Delta=\Delta_r$ of sizes $r\le r_0$ for some $r_0>0$. This is due to the fact that in the interior of the domain the solution is automatically in the class $W^{1,2}_{loc}(\Omega)$ implying that the estimate \eqref{E:bmo} will also holds on balls of sizes $r\ge r_0$ but with a slightly larger constant.\vglue1mm

\noindent{\it Remark 5.} We only assume that the condition \eqref{E:bmo} holds for all continuous data $f$. However, we claim that this implies the the same estimate holds for all bounded Borel measurable functions $f$ as a consequence. To see this it is enough to realize that if $f_j\to f$ in the sense that
$$\int f_j\,d\mu\to \int f\,d\mu,\quad\mbox{for any Borel probability measure }\mu\mbox{ on }\partial\Omega,$$
then if $u_j$ (or $u$) is the solution of the parabolic boundary value problem with data $f_j$ (or $f$), respectively, then for any compact set $K\subset \Omega$
we have $u_j\to u$ uniformly on $K$ as $j\to\infty$. Hence for any $\delta>0$ by Lemma \ref{L:Caccio} we have
$$\sigma^{-1}(\Delta_r)\int_{T(\Delta_r)\cap\{(X,t)\in\Omega:\delta(X,t)>\delta\}}|\nabla (u-u_j)|^2\delta\,dX \to 0,$$
and therefore
$$\sup_{\Delta\subset\partial\Omega}\sigma(\Delta)^{-1}\int_{T(\Delta)\cap\{(X,t)\in\Omega:\delta(X,t)>\delta\}}|\nabla u|^2\delta\,dX \lesssim \|f\|^2_{BMO(\partial \Omega, d\sigma)},$$
provided we have \eqref{E:bmo} for $u_j$ and $f_j$.
As this holds uniformly for all $\delta>0$ taking a limit $\delta\to 0$ yields \eqref{E:bmo} for $u$ and $f$. In particular, this implies that \eqref{E:bmo} holds for all bounded Borel measurable boundary data $f$.






\section{Basic results and Interior estimates}\label{S13:BI}

\smallskip





We now recall some estimates and tools needed for the proofs of Theorems \eqref{T:Main} and \eqref{T:Main2}.
\begin{lemma}\label{L:Caccio}(A Cacciopoli inequality, see \cite{Aron})
Suppose that $u$ is a weak solution of \eqref{E:v} in $\Omega$. For an interior point $(X, t) \in \Omega$  and any $0 <r < \delta(X,t)/4$ such that $Q_{4r}(X, t):=\{(Y,s):\,|X-Y|<4r\mbox{ and }|t-s|<16r^2\} \subset \Omega$, there exists a constant $C$ such that
\begin{equation*}\begin{split}
& r^{n} \left(\sup_{Q_{r}(X, t)} u \right)^{2} \\
&\leq    C \sup_{t-(2r)^2 \leq s \leq t+(2r)^2} \int_{B_{2r}(X)} u^{2}(Y,s) \,dY
       + C\int_{Q_{2r}(X, t)} |\nabla u|^{2} \,dY\,ds \\
&\leq    \frac{C^2}{r^2} \int_{Q_{4r}(X, t)} u^{2}(Y, s) \,dY\,ds.
\end{split}\end{equation*}
\end{lemma}


Lemmas 3.4 and 3.5 in \cite{HL} give us the following estimates for a weak solution of \eqref{E:v}.

\begin{lemma}\label{Holder}(Interior H\"{o}lder continuity)
Suppose that $u$ is a weak solution of \eqref{E:v} in $\Omega$. If $|u| \leq K < \infty$ for some constant $K>0$ in $Q_{4r}(X, t) \subset \Omega$, then for any $(Y, s), (Z, \tau) \in Q_{2r}(X, t)$ there exists a constant $C>0$ and $0 < \alpha < 1$ such that
\[
\left|u(Y, s) - u(Z, \tau)\right| \leq C K \left( \frac{|Y-Z| + |s - \tau|^{1/2}}{r} \right)^{\alpha}.
\]
\end{lemma}
\begin{lemma}
  \label{lemma_boundary-hoelder-continuity}(Boundary H\"older Continuity). Let
  $u$ be a weak solution of \eqref{E:v} in $T ( \Delta_{2 r} ( y, s))$. If $r > 0$
  and $u$ vanishes continuously on $\Delta_{2 r} ( y, s)$, then there exists
  $C$ and $\alpha$, $0 < \alpha < 1 \leqslant C < \infty$, such that for $( X,
  t) \in T ( \Delta_{r / 2} ( y, s))$,
  \begin{eqnarray*}
    u ( X, t) = u ( x_0, x, t) \leqslant C ( x_0 / r)^{\alpha} \max_{T (
    \Delta_r ( y, s))} u.
  \end{eqnarray*}
  If $u \geqslant 0$ in $T ( \Delta_{2 r} ( y, s))$ then there exists $C$ such
  that for $( X, t) \in T ( \Delta_{r / 2} ( y, s))$,
  \begin{eqnarray*}
    u ( X, t) \leqslant C ( x_0 / r)^{\alpha} u ( r, y, s + 2 r^2) .
  \end{eqnarray*}
\end{lemma}

\begin{lemma}\label{Harnack}(Harnack inequality)
Suppose that $u$ is a weak nonnegative solution of \eqref{E:v} in $U$ such that $Q_{4r}(X, t) \subset U$. Suppose that $(Y,s), (Z, \tau) \in Q_{2r}(X,t)$. There exists an a priori constant $c$ such that, for $\tau < s$,
\[
u(Z, \tau) \leq u(Y, s) \exp \left[ c\left( \frac{|Y-Z|^2}{|s-\tau|} + 1 \right) \right].
\]
If $u\geq 0$ is a weak solution of the adjoint operator of \eqref{E:v}, then this inequality is valid when $\tau > s$.
\end{lemma}

We state a version of the maximum principle from \cite{DH}, that is a modification of Lemma 3.38 from \cite{HL}.
\begin{lemma}\label{L:MP}(Maximum Principle) Let $u$, $v$ be bounded continuous weak solutions to \eqref{E:v} in $\Omega$. If $|u|,|v|\to 0$ uniformly as $t\to-\infty$ and
\[
\limsup_{(Y,s)\to (X,t)} (u-v)(Y,s) \leq 0
\]
for all  $(X,t)\in\partial\Omega$, then $u\leq v$ in $\Omega$.
\end{lemma}

\noindent{\it Remark.} The proof of Lemma 3.38 from \cite{HL} works given the assumption that $|u|,|v|\to 0$ uniformly as $t\to-\infty$.  Even with this additional assumption,  the lemma as stated is sufficient for our purposes. We shall mostly use it when $u\le v$ on the boundary of $\Omega\cap\{t\ge \tau\}$ for a given time $\tau$. Obviously then the assumption that $|u|,|v|\to 0$ uniformly as $t\to-\infty$ is not necessary. Another case when the Lemma as stated here applies is when $u|_{\partial\Omega},v|_{\partial\Omega}\in C_0(\partial\Omega)$, where $C_0(\partial\Omega)$ denotes the class of continuous functions decaying to zero as $t\to\pm\infty$. This class is dense in any $L^p(\partial\Omega,d\sigma)$, $p<\infty$ allowing us to consider an extension of the solution operator from $C_0(\partial\Omega)$ to $L^p$.

\begin{lemma}
  \label{lemma_parabolicdoubling}(Parabolic doubling, corkscrew point, c.f. \cite{HL}, and \cite{N} for the doubling property of parabolic measure in time-varying domains). 
   Let $\Delta_{2 r} \subset \Delta_d$ be boundary balls  and let $V_{2r}$ and $V_d$
be their corkscrew points.
  Let
  $\omega^{V_d}$ be parabolic measure.   There exists $c > 1$ (depending only on the domain and the parabolic operator $\partial_t-\di(A\nabla\cdot)$) such that
  \begin{enumerate}
    \item[a)] $c \omega^{V_d} (\Delta_d) \ge 1$

    \item[b)] $\omega^{V_d} (\Delta_{2 r}) \le c \omega^{V_d}
    (\Delta_r)\qquad$ (doubling)

    \item[c)] If $E \subset \Delta_{2 r}$ is a Borel set and $\omega^{V_{2r}} (E) \ge \eta$, then $c \omega^{V_d} (E)
    \ge \eta \omega^{V_d} (\Delta_{2 r}) .$
  \end{enumerate}
\end{lemma}


\section{Proof of Theorem \ref{T:Main}}

We shall establish that the estimate
$$\sup_{\Delta\subset\partial\Omega}\sigma(\Delta)^{-1}\int_{T(\Delta)}|\nabla u|^2\delta\,dX \lesssim \|f\|^2_{BMO(\partial \Omega, d\sigma)}$$
holds for all solutions $u$ in $\Omega$ with boundary data $f$. As we have noted above it suffices to show this result for all balls $\Delta=\Delta_r$
of diameter $\le r'$ for some $r'>0$.

Consider any boundary ball $\Delta_r=\Delta_r(y,s)$ of size $\le r'$. Let $\Delta_{2^jr}(y,s)$ for $j\ge 0$ be the $2^j$-fold enlargement of the original ball $\Delta_r$. We want to consider all $j\le m$ where $m$ is the smallest integer such that
$$(X,t)\in\Omega\setminus \Delta_{2^mr}\quad\Longrightarrow\quad |t-s|\ge r_0^2,$$
where $r_0$ is the scale from Definition \ref{D:domain}. From now on we denote $d=2^mr$.

Let $f$ be a BMO function on $\partial\Omega$ and let $u$ be the unique solution of the boundary value problem with boundary data $f$.
We decompose $f$ into several pieces. First,
  decompose it into near and far parts:
  \begin{eqnarray*}
    f &=& ( f - \langle f \rangle_{\Delta_{2r}}) \chi_{\Delta_{2r}} +  ( f - \langle f \rangle_{\Delta_{2r}}) \chi_{\Delta_d\setminus\Delta_{2r}}
    +   ( f - \langle f
    \rangle_{\Delta_{2r}}) \chi_{\partial \Omega \backslash \Delta_{d}} + \langle f
    \rangle_{\Delta_{2r}} \\&=&f_1+f_2 + f_3 + \langle f
    \rangle_{\Delta_{2r}}.
  \end{eqnarray*}
Here and throughout the paper we use the notation $\langle f\rangle_B=\sigma(B)^{-1}\int_B f\,d\sigma$.   If $u_i$ is the solution for boundary data $f_i$, then we have for the solution $u$ with data $f$, $\nabla u = \sum^3_{i = 1} \nabla u_i$. We estimate the contribution of
  each $u_i$ separately. Observe that the term $\langle f
    \rangle_{\Delta_{2r}}$ plays no further role as the solution corresponding to it is constant and hence has zero gradient.\vglue1mm

We start with $u_1$, the solution for boundary data $f_1$. For a fixed point
  $( X, t) \in T ( \Delta_r)$, let us consider the set
  $\mathcal{T}^{\alpha}_{( X, t)} = \{ ( z, \tau) \in \Delta_r : ( X, t) \in
  \Gamma_{\alpha} ( z, \tau) \}$. Note that $\sigma ( \mathcal{T}^{\alpha}_{( X,
  t)}) \sim \delta ( X, t)^n$ with constant dependent on $\alpha$.
  \begin{eqnarray*}
    &  & \frac{1}{\sigma ( \Delta_r)} \int_{T ( \Delta_r)} | \nabla u_1 |^2
    \delta ( X, s) d X d s\\
    & \lesssim & \frac{1}{\sigma ( \Delta_r)} \int_{T ( \Delta_r)} | \nabla
    u_1 |^2 \delta ( X, s)^{- n + 1} \sigma ( \mathcal{T}_{( X, s)}) d X d s\\
    & \lesssim & \frac{1}{\sigma ( \Delta_r)} \int_{( x, s) \in \Delta_r}
    \int_{\Gamma_{\alpha, r} ( x, s)} \delta ( X, s)^{- n + 1} | \nabla u_1
    |^2 d X d s d \sigma\\
    & \lesssim & \frac{1}{\sigma ( \Delta_r)} \int_{\Delta_r} \left(S^r_{\alpha}\right)^2
    ( u_1) d \sigma.
\end{eqnarray*}
Since we assume that the parabolic measure
  to belong to $A_{\infty}$, there exists a (large) $p<\infty$ for which the Dirichlet problem is
solvable in $L^p$. Hence by H{\"o}lder's inequality and solvability (see the Remark 3 below Definition \ref{D:infty})
\begin{eqnarray*}
    \frac{1}{\sigma ( \Delta_r)} \int_{\Delta_r} \left(S^r_{\alpha}\right)^2 ( u_1) d
    \sigma & \lesssim & \left( \frac{1}{\sigma ( \Delta_r)} \int_{\Delta_r}
    \left(S^r_{\alpha}\right)^p ( u_1) d \sigma \right)^{2 / p}\\
    & \lesssim & \frac{1}{\sigma ( \Delta_r)^{2 / p}} \left( \int_{\Delta_{2r}}
    | f_1 |^p d \sigma \right)^{2 / p} .
\end{eqnarray*}

Since $f_1 = f - \langle f \rangle_{\Delta_1}$ on $\Delta_1$, this is the
  BMO estimate with exponent $p$. John-Nirenberg's inequality
$$\frac{1}{\sigma ( \Delta_{2r})^{1 / p}} \left( \int_{\Delta_1}
    | f-\langle f\rangle_{\Delta_{2r}} |^p d \sigma \right)^{1 / p}\lesssim \|f\|_{BMO}$$
and the fact that $\sigma$ is doubling gives us our desired estimate.

In order to estimate the contribution of $u_2$, we write $f_2 = f^+_2 -
f^-_2$ with $f^{\pm}_{^{} 2} \geqslant 0$. Denote by $u^{\pm}_2$ the
corresponding solutions. By linearity of the equation we have for the
solution $u_2$ with data $f_2$ that $u_2 = u^+_2 - u^-_2$, whose contributions we
estimate separately. Let us denote by $\tilde{u}$ the solution with boundary
data $|f_2|$.

We now cover $T (\Delta_r)$ by a union of balls $B_i$, $i\in\mathbb N$, of finite overlap with the following properties.
$$T(\Delta_r)\subset\bigcup_i B_i\subset \bigcup_i (1+\delta)B_i \subset T(\Delta_{3/2r}).$$
Here $(1+\delta)B_i$ is a small enlargement of the ball $B_i$. All points
$( X, t) \in (1+\delta)B_i$ have comparable distance to
  the boundary, specifically, $\delta ( X, t) \sim \mbox{diam}(B_i)$. Furthermore, each point $(X,t)\in \bigcup_i (1+\delta)B_i$
is covered by at most $K$ enlarged balls $(1+\delta)B_i$, where $K$ only depends on the character of the
admissible domain $\Omega$. We have
\begin{eqnarray}
    &  & \frac{1}{\sigma ( \Delta_r)} \int_{T ( \Delta_r)} | \nabla u^{\pm}_2
    |^2 \delta ( X, t) d X d t\nonumber\\
    & \le & \frac{1}{\sigma ( \Delta_r)} \sum_{i}
    \int_{B_i} | \nabla u^{\pm}_2 |^2 \delta ( X, t) d X dt \nonumber\\
    & \lesssim & \frac{1}{\sigma ( \Delta_r)} \sum_{i} ( \mbox{diam} (B_i)^{- 1} \int_{( 1 + \delta) B_i}
    {u_2^{\pm}}^2 ( X, t) d X d t\nonumber\\
    & \lesssim & \frac{1}{\sigma ( \Delta_r)} \sum_{i} \int_{( 1 + \delta) B_i} {u_2^{\pm}}^2 ( X, s) \delta ( X,
    t)^{- 1} d X d t\label{E:imp}\\
    & \lesssim & \frac{K}{\sigma ( \Delta_r)} \int_{T ( \Delta_{3/2r})}
    \tilde{u}^2 ( X, t) \delta ( X, t)^{- 1} d X d t\nonumber\\
    & \lesssim & \| f \|^2_{{BMO}} \frac{r^{- 2 \varepsilon}}{\sigma (
    \Delta_r)} \int_{T ( \Delta_{3/2r})} \delta ( X, t)^{2 \varepsilon - 1} d X d t\nonumber
  \end{eqnarray}
  We use Lemma \ref{L:Caccio}, a
  pointwise estimate $u^{\pm}_2 \leqslant \tilde{u}$ as well as the pointwise
  estimate of  Lemma \ref{lemma_pointwise-BMO-estimates} below for the solution $\tilde{u}$ in terms of $\| f \|_{{BMO}}$ which uses the fact that $\tilde{u}$ vanishes on the boundary:
  
Observe now that the last
  expression is summable in the sense that
  \begin{eqnarray*}
    r^{- 2 \varepsilon} \int_{T ( \Delta_{3/2r})} \delta ( X, t)^{2 \varepsilon -
    1} d X d t \sim \sigma ( \Delta_r) .
  \end{eqnarray*}
  Therefore we get  a bound of  $\| f \|_{{BMO}}^2$ for the term with $u_2$.  
 
  We now state a lemma whose proof we postpone to the end.

\begin{lemma}
  \label{lemma_pointwise-BMO-estimates} Let $\Delta_d$ be a boundary cube of
  scale comparable to the diameter of the domain and $\Delta_{4 r} \subset
  \Delta_d$. Let $u$ denote the solution to boundary data $| f |
  \chi_{\Delta_d\setminus\Delta_{2r}}$, where $f$ is a function with $\langle f \rangle_{\Delta_{2r}}=0$. There exists $\varepsilon >
  0$ such that $\forall ( X, t) \in T ( \Delta_r)$,
  \begin{eqnarray*}
    u ( X, t) \lesssim \| f \|_{{BMO}} \qquad\mbox{and}\qquad u ( X, t)
    \lesssim \left( \tfrac{\delta ( X, t)}{r} \right)^{\varepsilon} \| f
    \|_{{BMO}} .
  \end{eqnarray*}
\end{lemma}

For the solution $u_3$ with boundary data $f_3$ we consider a further decomposition. 
(Up to this point, the argument has been exactly as in the elliptic case (\cite{DKP}).)
For all $j\ge 1$ consider
$$U_j=\{(X,t)\in\Omega\setminus \Delta_d:\,t\in [s-jr_0^2,s-(j+1)r_0^2) \}.$$
Since the scale $r_0$ is comparable to the diameter of each time slice of $\Omega$ each $U_j$ is contained in some boundary ball $\Delta^j$
of radius comparable to $r_0$ (and $d$) with $\sigma(U_j)\approx \sigma(\Delta^j)\approx r_0^n$.

We write
$$f_3=\sum_{j\ge 1}(f-\langle f\,\rangle_{\Delta_{2r}})\chi_{U_j}+h=\sum_{j\ge 1} g_j+h.$$
here $h$ is the portion of $f_3$ supported on $\Omega\cap\{t-s\ge r_0^2\}$. Observe that the term $h$ plays no further role as we only need to prove the estimate for $u$ on $T(\Delta_r)$, where the contribution from $h$ is zero. Hence it remains to deal with the data $g_j$, we denote the corresponding solutions by $w_j$.
We estimate the $L^p$ norm of $w_j$. We can add and subtract constants in order to use the BMO condition, writing $g_j$ on $U_j$ as:
$$g_j=(f-\langle f\,\rangle_{\Delta^{j}}) + \sum_{k=2}^j (\langle f\rangle_{\Delta^{k}}-\langle f\rangle_{\Delta^{k-1}})+(\langle f\rangle_{\Delta^{1}}-\langle f\rangle_{\Delta_{d}})+\sum_{k=2}^m(\langle f\rangle_{\Delta_{2^kr}}-\langle f\rangle_{\Delta_{2^{k-1}r}}).$$
The BMO condition on $f$ entails that for a ball of any radius $s$, and its double:
$$\|\langle f\rangle_{\Delta_{2s}}-\langle f\rangle_{\Delta_{s}}\|_{L^\infty}\lesssim\|f\|_{BMO},$$
and hence
$$\|g_j-(f-\langle f\,\rangle_{\Delta^{j}})\|_{L^\infty}\lesssim (j+m)\|f\|_{BMO}.$$
Again by John-Nirenberg we have for $f-\langle f\,\rangle_{\Delta^{j}}$ and any $p>1$ on $\Delta^j$:
$$\frac{1}{\sigma ( \Delta^{j})^{1 / p}} \left( \int_{\Delta^j}
    | f-\langle f\rangle_{\Delta^{j}} |^p d \sigma \right)^{1 / p}\lesssim \|f\|_{BMO},$$
and therefore for any $p>1$ we have that
$$\|g_j\|_{L^p(U_i)}\lesssim \sigma ( \Delta^{j})^{1 / p}(j+m+1)\|f\|_{BMO}.$$
The $A_\infty$ assumption as already noted above implies $L^p$ solvability of the Dirichlet boundary value problem for some large $p$; in particular this gives
that $\|N(w_j)\|_{L^p(\partial\Omega)}\lesssim \sigma( \Delta^{j})^{1 / p}(j+m+1)\|f\|_{BMO}$. Recall that $g_j=0$ for all times larger than $s-jr_0^2$. 

Since $w_j$ vanishes on the boundary of $U_{j-1}$, boundary H\"older regularity gives 
$$ \|w_j\|_{\infty} \lesssim \frac{1}{r_0 \sigma ( \Delta^{j})^{1 / p}} \left( \int_{U_{j-1}}
    | w_j |^p dXdt \right)^{1 / p}.$$

The solid integral over $U_{j-1}$ can be dominated by a nontangential maximal function:

$$ \frac{1}{r_0 \sigma ( \Delta^{j})^{1 / p}} \left( \int_{U_{j-1}}
    | w_j |^p dXdt \right)^{1 / p} \lesssim  \frac{1}{\sigma ( \Delta^{j})^{1 / p}} \left( \int_{\Delta^j}
    | N(w_j) |^p d \sigma \right)^{1 / p}.$$

And the above estimate for the nontangential maximal function of $w_j$ implies 
$$\|w_j\|_{L^\infty\left(\Omega_{s-(j-1)r_0^2}\right)}\lesssim (j+m+1)\|f\|_{BMO}.$$
Here, as before, $\Omega_\tau$ denotes the time slice at time $\tau$ of the domain $\Omega$. 

Now we are able to use the exponential decay for the solution
of any parabolic PDE with vanishing Dirichlet data on the lateral boundary. It follows that
$$\|w_j\|_{L^\infty(\Omega_t)}\lesssim e^{\beta(s-(j-1)r_0^2-t)}(j+m+1)\|f\|_{BMO},\qquad\mbox{for all }t\ge s-(j-2)r_0^2,$$
where the decay parameter $\beta>0$ only depends on the ellipticity constants and $\sup_{t}\diam(\Omega_t)$. In particular for $T(\Delta_d)$ we get
$$\|w_j\|_{L^\infty(T(\Delta_d))}\lesssim e^{-\beta jr_0^2}(j+m+1)\|f\|_{BMO}.$$

Finally, we use this to get an $L^\infty$ estimate on $T(\Delta_{2r})$. From Lemma \ref{lemma_boundary-hoelder-continuity} on the ball $\Delta_d$ with $x_0\le 2t$ we obtain
$$\|w_j\|_{L^\infty(T(\Delta_{2r}))}\lesssim 2^{-\alpha m}e^{-\beta jr_0^2}(j+m+1)\|f\|_{BMO}.$$
This final estimate allow us to do the same calculation as \eqref{E:imp} for $u_2^\pm$. We obtain:
$$\frac{1}{\sigma ( \Delta_r)} \int_{T ( \Delta_r)} | \nabla w^{\pm}_j
    |^2 \delta ( X, t) d X d t \lesssim 2^{-2\alpha m}e^{-2\beta jr_0^2}(j+m+1)\|f\|^2_{BMO}.$$
Finally, as $2^{-2\alpha m}m$ can be bounded independent of $m$ and $e^{-2\beta jr_0^2}j$ can be summed over all $j\ge 1$ we get for $u_3$:
\begin{eqnarray*}
&&\frac{1}{\sigma ( \Delta_r)} \int_{T ( \Delta_r)} | \nabla u_3
    |^2 \delta ( X, t) d X d t \\
&\lesssim& \sum_{j\ge 1} \frac{1}{\sigma ( \Delta_r)} \int_{T ( \Delta_r)} [|\nabla w_j^+|^2+|\nabla w_j^-|^2] \delta ( X, t) d X d t\\
&\lesssim & \|f\|_{BMO}^2\left(\sum_{j\ge 1} e^{-2\beta jr_0^2}j\right)\lesssim \|f\|_{BMO}^2.
\end{eqnarray*}
This concludes the proof of Theorem \ref{T:Main}, apart from the proof of Lemma \ref{lemma_pointwise-BMO-estimates}.\qed
\vglue2mm

\noindent{\it Proof of Lemma \ref{lemma_pointwise-BMO-estimates}.} We first note that the estimate
$$u ( X, t)
    \lesssim \left( \tfrac{\delta ( X, t)}{r} \right)^{\varepsilon} \| f
    \|_{{BMO}}$$
follows from the estimate $u ( X, t)\lesssim\|f\|_{{BMO}}$ by applying Lemma \ref{lemma_boundary-hoelder-continuity}, hence we shall only establish this bound.

We fix a corkscrew point of the ball $\Delta_d$ and denote it by $V_d$. Recall that
$\langle f\rangle_B$ denotes the average of $f$ over a ball $B$ with respect to the surface measure $\sigma$. As we want to consider averages with respect to
the parabolic measure $\omega$ as well, we use the notation
$$\langle f\rangle_{\omega,B}=\omega(B)^{-1}\int_B f\,d\omega,$$
occasionally using $\langle f\rangle_{\omega^{V_d},B}$ as well if we have to emphasize with respect to which point $\omega^{V_d}$ is defined.

We would like to replace the assumption $\langle f\rangle_{\Delta_{2r}}=0$  by $\langle f\rangle_{\omega^{V_d},\Delta_{2r}}=0$. We can do that by considering
a boundary value problem with data $|g|{\chi_{\Delta_d\setminus\Delta_{2r}}}$  instead, where
$g=f-\langle f\rangle_{\omega^{V_d},\Delta_{2r}}$. Denote the solution of this boundary value problem by $v$. If follows by the maximum principle that
$$\|u-v\|_{L^\infty(\Omega)}\le \|f-g\|_{L^\infty(\partial\Omega)}\le \left|\langle f\rangle_{\omega^{V_d},\Delta_{2r}}\right|.$$
Because the measures $\sigma$ and $\omega^{V_d}$ are $A_\infty$ with respect to each other, arguments exactly as in \cite{FaNeri}  entail that the difference of averages $\langle f\rangle_{\Delta_{2r}}$
and $\langle f\rangle_{\omega^{V_d},\Delta_{2r}}$ satisfy the following:
$$\left|\langle f\rangle_{\omega^{V_d},\Delta_{2r}}\right|=\left|\langle f\rangle_{\omega^{V_d},\Delta_{2r}}- \langle f\rangle_{\Delta_{2r}}\right|\lesssim \|f\|_{BMO},$$
This gives $\|u-v\|_{L^\infty(\Omega)}\lesssim\|f\|_{BMO}$, hence if $v ( X, t)\lesssim\|g\|_{{BMO}}=\|f\|_{{BMO}}$ for $(X,t)\in T(\Delta_r)$ then we have $u ( X, t)\lesssim\|f\|_{{BMO}}$ as well. From now of we can therefore assume that $\langle f\rangle_{\omega^{V_d},\Delta_{2r}}=0$.

For $j\ge 2$ we consider dyadic annuli $S_j=\Delta_{2^jr}\setminus \Delta_{2^{j-1}r}$.

\begin{eqnarray*}
    u ( Z, \tau) & \leqslant & \sum_{j=2}^m \int_{S_j} \left| f -
    \langle f \rangle_{\omega^{V_d}, \Delta_{2^jr}} \right| K^{( Z, \tau)} ( y, s)
    d y d s\\
    &  & + \sum_{j=2}^m \left| \langle f \rangle_{\omega^{V_d},
    \Delta_{2^jr}} \right| \int_{S_j} K^{( Z, \tau)} ( y, s) d y d s.
\end{eqnarray*}

Here $K^{( Z, \tau)}$ denotes the Radon-Nikodyn derivative of the
  parabolic measure at point $( Z, \tau)$, i.e., $d \omega^{( Z, \tau)} ( y, s) = K^{( Z, \tau)} ( y, s)d\sigma(y,z)$.

We have for the first sum:
  \begin{eqnarray*}
    &  & \sum_{j \geqslant 2} \int_{S_j} \left| f - \langle f
    \rangle_{\omega^{V_d}, \Delta_{2^jr}} \right| K^{( Z, \tau)} ( y, s) d y d s\\
    & \leqslant & \sum_{j \geqslant 2} \int_{\Delta_{2^jr}} \left| f - \langle f
    \rangle_{\omega^{V_d}, \Delta_{2^jr}} \right| K^{( Z, \tau)} ( y, s)\,dy ds\\
    & \lesssim & \sum_{j \geqslant 2} 2^{- \alpha j} \frac{1}{\omega^{V_d} ( \Delta_{2^jr})} \int_{\Delta_{2^jr}} \left| f - \langle f
    \rangle_{\omega^{V_d}, \Delta_{2^jr}} \right| d \omega^{V_d} ( y, s)\\
    & \lesssim & \| f \|_{{BMO}, \omega^{_{V_d}}} \sum_{j \geqslant 2}
    2^{- \alpha j} \lesssim  \| f \|_{{BMO}} .
  \end{eqnarray*}
  In the second line, we use Lemma \ref{lemma_change-variable-kernel} (see
  below). In the last line we use equivalence of BMO norms with respect to
  parabolic and surface measures that holds due to the $A_\infty$ assumption we have made. Finally, $V_r$ denotes the corkscrew point of the ball $\Delta_r$.

For the second sum we have (using $\langle f \rangle_{\omega^{V_d},
    \Delta_{2r}}=0$):
  \begin{eqnarray*}
    &  & \sum_{j=2}^m \left| \langle f \rangle_{\omega^{V_d},
    \Delta_{2^jr}} \right| \int_{S_j} K^{( Z, \tau)} ( y, s) d y d s
     =  \sum_{j=2}^m \left| \langle f \rangle_{\omega^{V_d},
    \Delta_{2^jr}} \right| \int_{S_j} d \omega^{( Z, \tau)} ( y, s)\\
    & \leqslant & \sum_{j \geqslant 2} \left( \sum_{i = 2}^j \left| \langle f
    \rangle_{\omega^{V_d}, \Delta_{2^ir}} - \langle f \rangle_{\omega^{V_d},
    \Delta_{2^{i-1}r}} \right| \right) \int_{S_j} d \omega^{( Z, \tau)} ( y, s)\\
    & \leqslant & \sum_{j=2}^m j \sup_i \left| \langle f
    \rangle_{\omega^{V_d}, \Delta_{2^ir}} - \langle f \rangle_{\omega^{V_d},
    \Delta_{2^{i-1}r}} \right| \int_{S_j} d \omega^{( Z, \tau)} ( y, s)\\
    & \lesssim & \| f \|_{{BMO}, \omega^{V_d}} \sum_{j \geqslant 0} j
    \int_{S_j} d \omega^{( Z, \tau)} ( y, s)\\
    & \lesssim & \| f \|_{{BMO}, \omega^{V_d}} \sum_{j \geqslant 0} j
    2^{- \alpha j} \int_{\Delta_{2^jr}} \frac{1}{\omega^{V_d} ( \Delta_{2^jr})} d \omega^{V_d} ( y, s)\\
    & = & \| f \|_{{BMO}, \omega^{V_d}} \sum_{j=2}^m j 2^{-
    \alpha j} \lesssim \| f \|_{{BMO}}.
  \end{eqnarray*}
Here we again have used Lemma \ref{lemma_change-variable-kernel} stated below. Hence the result holds.\qed

\begin{lemma}
  \label{lemma_change-variable-kernel} Let, $d\lesssim r_0$, $\Delta_d \subset \partial \Omega$
  and $V_d$ the corkscrew point of $\Delta_d$. Let $\Delta_r \subset
  \Delta_d$ and denote by $V_r$ the cork screw point of $\Delta_r$. We have
  for each $j\ge 1$ such that $\Delta_{2^jr}\subset\Delta_d$
  \begin{eqnarray*}
    \sup_{( y, s) \in \Delta_{2^jr}} \sup_{( Z, \tau)
    \in T ( \Delta_r)} \frac{K^{( Z, \tau)} ( y, s)}{K^{V_d} ( y, s)} \lesssim
    2^{- \alpha j} \frac{1}{\omega^{V_d} ( \Delta_{2^jr})} .
  \end{eqnarray*}
\end{lemma}

\noindent {\it Proof.}   First observe that from assertion c) in Lemma
  \ref{lemma_parabolicdoubling} follows that for $E \subset \Delta_r$ we
  have
  \begin{eqnarray*}
    \omega^{V_r} ( E) \lesssim \frac{\omega^{V_d} ( E)}{\omega^{V_d} ( \Delta_r)} .
  \end{eqnarray*}
  Let $S_j$ be as before the dyadic annuli $\Delta_{2^jr} \backslash
  \Delta_{2^{j-1}r} .$ Let $V_j$ be the corkscrew point of the surface cube
  $\Delta_{2^jr}$. We apply the above inequality for an infinitesimally small cube (that is $t$ is tiny)
  $\Delta_t \subset S_j \subset \Delta^j$. We get
  \begin{eqnarray*}
    \omega^{V_j} ( \Delta_t) \lesssim \frac{\omega^{V_d} ( \Delta_t)}{\omega^{V_d} (
    \Delta_{2^jr})}\quad\mbox{and thus}\quad \frac{\omega^{A_j} ( \Delta_t)}{\omega^{V_d} ( \Delta_t)} \lesssim \frac{1}{\omega^{V_d}
    ( \Delta_{2^jr})} .
      \end{eqnarray*}
Now for $( Z, \tau) \in T ( \Delta_r)$ and $( y, s) \in \Delta_t \subset
  S_j$
  \begin{eqnarray*}
    \frac{K^{( Z, \tau)} ( y, s)}{K^{V_d} ( y, s)} = \lim_{t \rightarrow 0}
    \frac{\omega^{( Z, \tau)} ( \Delta_t)}{\omega^{V_d} ( \Delta_t)} = \lim_{t
    \rightarrow 0} \frac{\omega^{( Z, \tau)} ( \Delta_t)}{\omega^{A_j} ( \Delta_t)}
    \cdot \frac{\omega^{A_j} ( \Delta_t)}{\omega^{V_d} ( \Delta_t)} .
  \end{eqnarray*}
  Note that $\mbox{dist} ( \Delta_t, \Delta_r) \sim 2^jr$. Applying Boundary
  H{\"o}lder, Lemma \ref{lemma_boundary-hoelder-continuity}, using that the
  boundary data for $\omega^{( Z, \tau)} ( \Delta_t)$ vanishes on $\Delta_{2^{j-1}r}$ we have
  \begin{eqnarray*}
    \omega^{( Z, \tau)} ( \Delta_t) \lesssim \left( \frac{\mbox{dist} ( ( Z, \tau),
    \partial\Omega)}{2^jr} \right)^{\alpha} \omega^{A_j} ( \Delta_t) \lesssim 2^{- j
    \alpha} \omega^{A_j} ( \Delta_t)
  \end{eqnarray*}
  and the lemma is completely proved. \qed


\section{Proof of Theorems \ref{T:Main2} and \ref{T:Main3}}

We focus primarily on the parabolic case, since the elliptic case is less complicated. We start by recalling
the existence of dyadic grid that can be constructed for any doubling measure (\cite{Chr}).\vglue1mm


Let $\Delta_d \subset \partial \Omega$ where $d$ is of the size comparable to the scale $r_0$ from Definition \ref{D:domain}.
As before let $V^d$ be the corkscrew point of $\Delta_d$. By Lemma \ref{lemma_parabolicdoubling}, the parabolic
measure $\omega^{V_d} $ has the doubling property in $\Delta_d$, therefore the
metric space $\Delta_d$ has a dyadic grid with the following properties.
$\mathcal{D} ( \Delta_d)
= \{ I_{j}^l : j \in \mathbb{Z}, l \in \mathcal{I}_j \}$ with $I_{j}^l \subset
\Delta_d$ and $\mathcal{I}_j$ an index set.  This dyadic grid possesses the following properties:
\smallskip

\begin{enumerate}
  \item $\bigcup_l I_{j}^l = \Delta_d$; $\omega ( \partial I_{j}^l) = 0$ for all  $j, l$.

  \item $\emptyset \in \mathcal{D} ( \Delta_d)$; $\Delta_d \in \mathcal{D} (
  \Delta_d)$.

  \item $\mbox{int}(I_{j}^l) \cap \mbox{int}(I_{j}^{l'}) = \emptyset$ if $l \neq l'$. (Here int$(B)$ means the interior of the set $B$).

  \item There exist $( x_l, t_l)$, called the center of $I_{j}^l$, so that
  $\Delta_{2^{- j}} ( x_l, t_l) \subseteq I_{j}^l \subseteq \Delta_{M 2^{-
  j}} ( x_l, t_l)$ where $M$ only depends on doubling constant of $\omega$. 
  \item If $j \geq j'$ then  $I_{j}^l
  \subseteq I_{j'}^{l'}$ or $I_{j}^l
  \cap I_{j'}^{l'} = \emptyset$.

  \item When $I_{j}^l
  \subsetneq I_{j'}^{l'}$ then there exists $C < 1$ 
  so that $\omega ( I_{j}^l) < C \omega ( I_{j'}^{l'})$.

  \item Any open set $\mathcal{O} \subset \Delta_d$ can be decomposed as
  $\mathcal{O}= \bigcup_{j, l} I_{j}^l $ where $\mbox{int}(I_{j}^l)$ are pairwise
  disjoint and for each $I_j^l$, there is a point  $P_{j}^l \in \Delta_d \backslash
  \mathcal{O}$ such that $\mbox{dist} ( P_{j}^l, I_{j}^l) \sim \mbox{diam} (
  I_{j}^l) \sim 2^{- j}$.
\end{enumerate}

{\it Remark. If $S$ is an element of the dyadic grid we shall say that $S$ has scale $j$ if $S = I_j^l$ for some $l$.}

\begin{definition}
  \label{definition epsiloncover} (c.f. \cite{KKPT})  Let $\varepsilon_0$ be given. Let $E \subset
  \Delta_r \subseteq \Delta_d$. A good $\varepsilon_0$-cover for $E$ in
  $\Delta_r$ of length $k$ is a collection of nested open sets $\{
  \mathcal{O}_l \}^k_{l = 1}$ with $E \subseteq \mathcal{O}_k \subseteq \ldots
  \subseteq \mathcal{O}_0 = \Delta_r$. Moreover, each $\mathcal{O}_l$ decomposed as $\mathcal{O}_l=
  \bigcup^{\infty}_{i = 1} S^l_i$ such that
  \begin{enumerate}
    \item $S^l_i \in \mathcal{D} ( \Delta_d)\quad \forall i, l$

    \item $\omega^{V_d} ( \mathcal{O}_l \cap S^{l - 1}_i) \leqslant
    \varepsilon_0 \omega^{V_d} ( S^{l - 1}_i)\quad \forall\, 1 \leqslant l \leqslant
    k$.
  \end{enumerate}
\end{definition}

Note  that when $k \geqslant l > m > 0$, then $\omega ( S^m_j \cap \mathcal{O}_l)
\leqslant \varepsilon_0^{l - m} \omega ( S^m_j)$.

\begin{lemma}(c.f. Lemma 2.6 of \cite{KKPT})
  Let $E \subset \Delta_r \subset \Delta_d$. Given $\varepsilon_0 > 0$, there
  exists $\delta_0 > 0$ such that if
  $$\omega^{V_d} ( E) / \omega^{V_d} (
  \Delta_r) \leqslant \delta_0$$
  then $E$ has good $\varepsilon_0$-cover
  of length $k = k ( \varepsilon_0, \delta_0)$. In fact, $k \sim
  -\varepsilon_0 \log \delta_0$.
\end{lemma}

As explained in Remark 5, we may assume that for all Borel-measurable bounded $f$, the solution $u$ with boundary data $f$ satisfies:
  \begin{eqnarray*}
    \sigma ( \Delta_r)^{- 1} \int_{T ( \Delta_r)} | \nabla u ( Y, s) |^2
    \delta ( Y, s) d Y d s \lesssim \| f \|_{\infty},
  \end{eqnarray*}
uniformly for all balls $\Delta_r\subset \partial\Omega$ with $r\le r'$ for some $r'>0$. As we have noted above this condition is equivalent to
saying that the truncated square function $S^r$ satisfies
  \begin{eqnarray*}
    \int_{\Delta_r} ( S^r(u))^2 d\sigma \leqslant \sigma ( \Delta_r),
  \end{eqnarray*}
for all $\|f\|_{L^\infty}\le 1$ and $0<r\le r'$. Recall that our goal is to prove that for all $E\subset\Delta_r\subset\Delta_d$:
  \begin{eqnarray*}
    \omega^{V_d} ( E) / \omega^{V_d} ( \Delta_r) < \delta \quad\Longrightarrow\quad \sigma (
    E) / \sigma ( \Delta_r) < \varepsilon .
  \end{eqnarray*}
We pursue the following strategy, as in \cite{KKKT}. We will establish that, given $\delta>0$, one can find $K(\delta)$ (with $K(\delta)\to \infty$ as $\delta\to0+$) such that
for some $f$ with $\|f\|_{L^\infty}\le 1$ we have for the solution $u$ corresponding to the boundary data $f$:
$$(S^r(u))^2(x,t)\ge K,\qquad\mbox{for all }(x,t)\in E.$$
This would imply that
$$K(\delta)\sigma(E)\le \int_{E}(S^r(u))^2d\sigma\le \int_{\Delta_r}(S^r(u))^2d\sigma \lesssim \sigma ( \Delta_r),$$
and hence
\begin{eqnarray*}
    \omega^{V_d} ( E) / \omega^{V_d} ( \Delta_r) < \delta \quad\Longrightarrow\quad \sigma (
    E) / \sigma ( \Delta_r) \le \frac{C}{K(\delta)},
  \end{eqnarray*}
from which $A_\infty$ follows as we choose $\delta>0$ such that $K(\delta)>C/\varepsilon$.  \vglue1mm

It remains to construct $f$ with the stated properties.  Assume therefore that $E\subset\Delta_r$ is given and that $\omega^{V_d} ( E) / \omega^{V_d} ( \Delta_r) < \delta$ where $\delta>0$ will be determined later.
Without loss of generality we may assume that $d\le \frac{r_0}{C_0}$ (c.f. Definition \ref{D:domain}) and hence our ball $\Delta_d$ is contained in one $\ell$-cylinder $\mathbb Z$ in which the boundary $\partial\Omega$ is given as a graph of a function $\phi$. On such local coordinate system we can simplify the geometry through the use of
a pull-back transformation which transforms our PDE into a new parabolic PDE on a subset of $U={\mathbb R}^+\times {\mathbb R}^{n-1}
\times \mathbb R$. 
Let $U=\mathbb{R}^+\times \mathbb{R}^{n-1} \times \mathbb{R}$.
We will consider a mapping $\rho : U \to \Omega$  known as the Dalhberg-Kenig-Stein adapted distance mapping, which appears also in Ne\v{c}as in the elliptic setting. In the parabolic setting this was  studied in \cite{HL}, and has been extensively used in a variety of contexts including \cite{KP}, \cite{Rn2}, and \cite{DH}. The mapping is given in local coordinates as follows:
\begin{equation}
\rho (x_0, x, t) = (x_0 + P_{\gamma x_0}\psi(x,t), x, t).\label{mapping}
\end{equation}
Here, $P(x,t) \in C_{0}^{\infty}(Q_{1}(0,0))$ is a non-negative function, defined for $(x,t) \in \mathbb{R}^{n-1}\times\mathbb{R}$, and
\[
P_{\lambda} (x,t) \equiv \lambda^{-(n+1)} P\left( \frac{x}{\lambda}, \frac{t}{\lambda^2} \right)
\]
and
\[
P_{\lambda} \psi(x,t) \equiv \int_{\mathbb{R}^{n-1}\times \mathbb{R}} P_{\lambda} (x-y, t-s) \psi(y,s) \,dy\,ds.
\]
Then $\rho$ satisfies
\[
\lim_{(y_0, y, s) \to (0, x, t)} P_{\gamma y_0} \psi(y,s) = \psi(x, t)
\]
and extends continuously to $\rho:\overline{U}\to \overline{\Omega}$. As follows from the discussion above the usual
surface measure on $\partial U$ is comparable with the measure $\sigma$ defined by \eqref{E:sigma} on $\partial\Omega$.
By setting $v = u \circ \rho$ and $f^v = f \circ \rho$, one finds that the equation \eqref{E:v} transforms to a new equation satisfied by $v$:
\begin{equation}\label{E:u}\begin{cases}
v_t = \di (A^{\rho} \nabla v) + \boldsymbol{B^{\rho}}\cdot \nabla v & \text{in } U, \\
v   = f                & \text{on } \partial U
\end{cases}\end{equation}
where $A^{\rho} = [a^{\rho}_{ij}(X, t)]$, $\boldsymbol{B}^{\rho} = [b^{\rho}_i (X, t)]$ are $(n\times n)$ and $(1\times n)$ matrices.

Hence for this new equation one can think of the ball $\Delta_r=\Delta_r(x,t)$ as the set
$$\{(0,y,s)\in U: |y-x|<r\mbox{ and }|s-t|<r^2\}.$$
The simplicity of this geometry is the primary reason for introducing the adapted distance mapping.

Consider a good $\varepsilon_0$-cover for $E$ relative to $\Delta_r$ ($\varepsilon_0$ to be determined).
This gives rise to sets $\{ \mathcal{O}_m
  \}^k_{m = 0}$ \ and to $\omega^{V_d}$-dyadic cubes $S_i^m$ so that
  $\mathcal{O}_m = \bigcup_i S_i^m$. When $m = 0$ then \ $\mathcal{O}_0 = \Delta_r$ and $S^0 =
  \Delta_r$. Notice that there exist $\sigma$-dyadic cubes so that $\Delta^m_i
  \subset S_i^m \subset M \Delta^m_i$ where the scales of $S^m_i$ and
  $\Delta^m_i$ are comparable to, say, $r^m_i$. We make the following
  convention on notation. If $\Delta^m_i = \Delta_{r^m_i} ( y_i^m, s_i^m)$ is
  said boundary cube, we denote by ${\Delta^m_i}'$ the boundary cube of scale
  \begin{eqnarray*}
    {r^m_i}' = r_i^m /\sqrt{2}
  \end{eqnarray*}
  centered at
  \begin{eqnarray*}
    ( y^m_i, {s^m_i}') = ( y^m_i, s^m_i - ( r^m_i)^2 / 2) = ( y^m_i, s^m_i - (
    {r^m_i}')^2) .
  \end{eqnarray*}

The function $f$ will be a sum of functions $f_m$ that we now define. For $m$ even and $0\le m < k$, we set
  \begin{eqnarray*}
    f_m(0,y,s) = \begin{cases} 1,&\quad \mbox{for } (0,y,s)\in \bigcup_i {\Delta^m_i}'\subset{\mathcal{O}_{m}},\\
                               0,&\quad\mbox{elsewhere}.
                 \end{cases}
  \end{eqnarray*}
  For $m$ odd we set
  \begin{eqnarray*}
    f_{m + 1} = - f_m \mathcal{X}_{\mathcal{O}_{m + 1}}.
  \end{eqnarray*}
  Observe that
  \begin{eqnarray*}
    f = \sum_{m = 0}^k f_m
  \end{eqnarray*}
  is a non-negative Borel function with $0\le f\le 1$.

Let us denote by $u$ the solution corresponding to $f$, by $u_m$ we denote the solutions with boundary data $f_m$.
We will show that data $f_m$ for $m$ an even integer, generates oscillation of square function of $u$
on a large enough subset $A_{a, m} ( x, t) \subset
  \Gamma_{a} ( x, t)$. Moreover, sufficiently many of the sets $A_{a, m} ( x, t)$ will be disjoint for distinct $m$.

Let $m$ be even. Take any $( x, t) \in E \subset \mathcal{O}_m = \bigcup_i
  S^m_i$ and find $S^m \in \{ S^m_i \}$ that contains $( x, t)$. Recall that
  $f_m = 1$ on ${\Delta^m}'$ and $f_m = 0$ elsewhere in $S^m$.
Let  $r^m$ be the scale of the radius of $\Delta^m$. If $( y^m,
  s^m)$ is the center of $\Delta^m$, then 
  \begin{eqnarray*}
    V_m' := ( {r^m}', y^m, s^m + ( {r^m}')^2)
  \end{eqnarray*}
  is the cork screw point of ${\Delta^m}'$. 
  We later choose $a$ to ensure that $V'_m \in
  \Gamma_{a} ( x, t)$. 

Clearly, by Lemma \ref{lemma_parabolicdoubling} part a),
  \begin{eqnarray*}
    \omega^{V_m'}({\Delta^m}')\gtrsim 1.
  \end{eqnarray*}
Since $f\ge f_m+f_{m+1}$ by the maximum principle we have
  \begin{eqnarray*}
    u ( {V_m}') & = & \int_{\Delta_r} f ( y, s) K^{{V_m}'} ( y, s) d y d s\\
    &\geqslant& \int_{\Delta_r} ( f_m ( y, s) + f_{m + 1} ( y, s))
    K^{V_m'} ( y, s) d y d s \\
    &\geqslant& \int_{\Delta^m} ( f_m ( y, s) + f_{m + 1} ( y, s))
    K^{V_m'} ( y, s) d y d s\\
    &=&\omega^{V_m'}({\Delta^m}')-\omega^{V_m'}({\Delta^m}'\cap {\mathcal O}_{m+1})\\
    &=& O(1)-\omega^{V_m'}({\Delta^m}'\cap {\mathcal O}_{m+1}),
  \end{eqnarray*}
where in the fifth line we have used the definition of $f_m$ and $f_{m+1}$. By Lemma \ref{lemma_parabolicdoubling} part c) for $E
  = {\Delta^m}' \cap \mathcal{O}_{m + 1}$ we have
  \begin{eqnarray*}
    \omega^{V_m'}({\Delta^m}'\cap {\mathcal O}_{m+1}) \lesssim
    \frac{\omega^{V_d} ( {\Delta^m}' \cap \mathcal{O}_{m + 1})}{\omega^{V_d} (
    S^m)} \leq \frac{\omega^{V_d} ( S^m \cap \mathcal{O}_{m +
    1})}{\omega^{V_d} ( S^m)} < C \varepsilon_0 .
  \end{eqnarray*}
 where the $C$ in the last line depends on doubling constants, and is independent of $m$.
It follows that for $\varepsilon_0$ chosen sufficiently small one has $u ( {V_m}')\gtrsim 1$. By the
Harnack inequality, Lemma \ref{Harnack}, 
there exists a substantial set of points later in time where this inequality holds. Namely
  \begin{eqnarray*}
    u ( {r^m}', y, s) \gtrsim 1 \qquad\mbox{for all } (y, s) \in H^m
  \end{eqnarray*}
  where
  \begin{eqnarray*}
    H^m := \{ | y^m - y | \leqslant {r^m}'' \} \times \left\{ \tfrac{1}{2}
    ( {r^m}'')^2 \leqslant s - ( s^m + ( {r^m}')^2) \leqslant ( {r^m}'')^2 \right\},
  \end{eqnarray*}
  with ${r^m}'' = {r^m}' / 4$. Again we postpone the considerations that will guarantee that $H^m \subset \Gamma_{a} (
  x, t)$. In the elliptic case,
similarly by the Harnack inequality, there is a small ball around ${V_m}'$ where $u \gtrsim 1$ holds.\vglue1mm

We now produce another set of points closer to the boundary than $H^m$ where $u$ is small. This
  will give us an estimate on the oscillation of the square function of $u$. Let us consider the values of  $u ( \rho {r^m}', y, s)$, for small $\rho$
and for  $( y, s) \in H^m$ .

Observe that by the construction of $f$,
  \begin{eqnarray*}
    f_m + \chi_{\mathcal{O}_{m + 1} \cap \Delta^m} \geqslant f \quad \mbox{on } \Delta^m ,
  \end{eqnarray*}
  since $\sum^{m - 1}_{l = 0} f_l = 0$ on $\Delta^m$ ($m$ is even).

As $f_m = 0$ on $\Delta_{{r^m}'} ( y^m, s^m + ( {r^m}')^2)$ we see that if we set $g=\chi_{\partial U \setminus \Delta_{{r^m}'} ( y^m, s^m + ( {r^m}')^2)}$
then $g+\chi_{\mathcal{O}_{m + 1} \cap \Delta^m}\ge f$ on $\partial U$.

Let $v$ be the solution with boundary data $g$. As $0\le v\le 1$ on $U$ by Lemma \ref{lemma_boundary-hoelder-continuity} we have an estimate for some $0 < \beta < 1$,
  \begin{eqnarray*}
    v ( \rho {r^m}', y, s) \lesssim \rho^{\beta} \quad\mbox{when } ( y, s) \in H^m .
  \end{eqnarray*}
  It remains to control the contribution of $\chi_{\mathcal{O}_{m + 1} \cap
  \Delta^m}$. We want to control $\omega^{( \rho {r^m}', y, s)}(\mathcal{O}_{m + 1} \cap
  \Delta^m)$ but we first look at $\omega^{V_m}(\mathcal{O}_{m + 1} \cap
  \Delta^m)$, where $V_m$ is the corkscrew point of $\Delta^m$.
 By Lemma   \ref{lemma_parabolicdoubling} part c) we get

 \begin{eqnarray*}
    \omega^{V_m} ( \mathcal{O}_{m + 1} \cap \Delta^m) &
    \lesssim & \frac{\omega^{V_d} ( \mathcal{O}_{m + 1} \cap
    \Delta^m)}{\omega^{V_d} (\Delta^m)}\\
    &\le & C_M \frac{\omega^{V_d} ( \mathcal{O}_{m + 1} \cap
    \Delta^m)}{\omega^{V_d} (M\Delta^m)}\\
    &\le & C_M \frac{\omega^{V_d} ( \mathcal{O}_{m + 1} \cap
    S^m)}{\omega^{V_d} (S^m)} \le C_{M} \varepsilon_0.
  \end{eqnarray*}
The dependence on $M$ stems from $\Delta^m \subset S^m \subset M \Delta^m$,
  where $M$ is an absolute constant.

Our goal is to bound $\omega^{( \rho {r^m}', y, s)}(\mathcal{O}_{m + 1} \cap
  \Delta^m)$ by some $\eta_2$ small. Assume instead that $\omega^{( \rho {r^m}', y, s)}(\mathcal{O}_{m + 1} \cap
  \Delta^m)\ge \eta_2$. Then by Lemma \ref{Harnack} since the time coordinate of $V_m$ is larger than $s$, one can construct a Harnack chain $B_1,B_2,\dots,B_j$ consisting of interior balls
  such that $B_i\cap B_{i+1}\ne\emptyset$, $( \rho {r^m}', y, s)\in B_1$, $V_m\in B_j$ and $4B_i\subset U$ for all $i$. The minimal length $j$ of such chain
  depends on $\rho$. By repeated application of Lemma \ref{Harnack} on each $B_i$ it follows that
$$C_{M} \varepsilon_0\ge\omega^{V_m} ( \mathcal{O}_{m + 1} \cap \Delta^m)\ge C_\rho\omega^{( \rho {r^m}', y, s)}(\mathcal{O}_{m + 1} \cap
  \Delta^m)\ge C_\rho \eta_2,$$
where $C_\rho$ is a small positive constant depending on $\rho$. If $\varepsilon_0$ is chosen small enough such that $\varepsilon_0<\frac{C_\rho}{C_m}\eta_2$ this is a contradiction. Hence we must have $\omega^{( \rho {r^m}', y, s)}(\mathcal{O}_{m + 1} \cap
  \Delta^m)< \eta_2$ and hence from
  $$v ( \rho {r^m}', y, s)+ \omega^{( \rho {r^m}', y, s)}(\mathcal{O}_{m + 1} \cap
  \Delta^m) \ge u( \rho {r^m}', y, s),$$
it follow that
$$u( \rho {r^m}', y, s) \le \eta_2 +C\rho^\beta.$$
Let us clarify the order in which we choose the parameters. Firstly, on $H^m$ we have $u({r^m}',\cdot,\cdot)\gtrsim 1$, or more precisely we have $u\ge 1-\eta_1$.
Choose $\rho>0$ such that $C\rho^\beta\le \frac13(1-\eta_1)$ and pick $0<\eta_2=\frac13(1-\eta_1)$. Then we choose $\varepsilon_0$ so small such that
$\varepsilon_0<\frac{C_\rho}{C_M}\eta_2$. This yields
$$u( \rho {r^m}', y, s)\le \textstyle\frac23(1-\eta_1),\quad\mbox{and}\quad u(  {r^m}', y, s)\ge 1-\eta_1,$$
for $(y,s)\in H^m$ and hence
$$|u(  {r^m}', y, s)-u( \rho {r^m}', y, s)|\ge \textstyle\frac13(1-\eta_1)>0,\qquad\mbox{for all } (y,s)\in H^m.$$
This is the key estimate that will allow us to show that the square function of $u$ is large on the set:
$$   A_m = \{ ( y_0, y, s) : \rho {r^m}' < y_0 < {r^m}', ( y, s) \in H^m
    \}$$
We claim that for large enough aperture $a$, $A_m \in \Gamma_a^{ r} ( x, t)$ for our initial point $(x,t)\in E$.
The choice of $a$ will depend on the $\ell$ in the character $(\ell,N,C_0)$ of the domain $\Omega$ and the ellipticity constant of the matrix $A$
and on $\rho$: the construction
above ensures that
$$r^m\approx \diam(A_m)\approx \mbox{dist}(A_m,\partial U)\approx \mbox{dist}(A_m,(x,t)),$$
where the implied constants depend on $\rho$ and $r^m$ denotes, as before, the scale of radius of the ball $\Delta^m\ni (x,t)$.

It follow that for any $( y, s) \in H^m$
  \begin{eqnarray*}
    1 & \lesssim & | u ( {r^m}', y, s) - u ( \rho {r^m}', y, s) |^2\\
    & = & \left| \int^1_0 \nabla u ( \rho {r^m}' + t ( {r^m}' - \rho
    {r^m}'), y, s) ( {r^m}' - \rho {r^m}') d t \right|^2\\
    & \leqslant & ( 1 - \rho)^2 ( {r^m}')^2\int^1_0 | \nabla u ( \rho r_{( m)} + t ( {r^m}' - \rho
    {r^m}'), y, s) |^2  d t\\
    & \lesssim & ( {r^m}')^{n + 1} \int^1_0 | \nabla u ( \rho {r^m}' + t ( {r^m}' -
    \rho {r^m}'), y, s) |^2 ( {r^m}')^{- n + 1} d t\\
    & \lesssim & ( {r^m}')^{n + 1} \int_{\rho {r^m}'}^{ {r^m}'} | \nabla u ( y_0), y, s) |^2 ( {r^m}')^{- n} d y_0.
  \end{eqnarray*}
  Hence integrating both sides over $H_m$ and dividing by $\sigma(H_n)\approx ({r^m}')^{n + 1}$ will give us
$$1\lesssim \int_{A_m} | \nabla u ( y_0), y, s) |^2 y_0^{- n} d y_0\, dy\, ds.$$
This  is the contribution of each $A_{m}$ to the lower bound on the square
function $S_r^a(u)(x,t)$ for $A_m\subset \Gamma^a_r(x,t)$. However, not all $A_m$ ($m$ even) are necessary disjoint.
To ensure $(A_{m_j})$ are disjoint we take a subsequence $m_j$ of even integers such that
$$\rho  {r^{m_j}}' > {r^{m_{j + 1}}}'.$$
Note that  by property (6), and for any level $m$, we have obtained a sequence, $S_j^m$, of elements of the dyadic grid with the
property that $S_j^m$ is properly contained in the $S_j^{m+1}.$ If $S_j^m$ has scale $2^{-i}$ then $S_j^{m+1}$ must have scale
at most $2^{-i-1}$, since by property (5) any two dyadic cubes at the same scale that are not identical are disjoint.
Thus, by skipping a fixed finite number of levels, 
choosing $m_{j + 1}={m_j}+2k$ for some fixed $k\in \mathbb N$
we see that 
$r^{m_{j + 1}} \leq M2^{-2k} r^{m_j } < \rho r^{m_j}$, with $M$ from property (4) and $k$ chosen such that $M2^{-2k}\le\rho$.

The number of disjoint $A_{m_j}$ is proportional to the length of good $\varepsilon_0$-cover, i.e.,\newline
 $\varepsilon_0
  \log \left( \omega^{V_d} ( \Delta_r) / \omega^{_{V_d}} ( E) \right)$. This, for $(x,t)\in E$ we have
$$(S^r)^2(u)(x,t)\gtrsim \sum_j \int_{A_{m_j}} | \nabla u ( y_0), y, s) |^2 y_0^{- n} d y_0\, dy\, ds\gtrsim \varepsilon_0
  \log \left( \omega^{V_d} ( \Delta_r) / \omega^{_{V_d}} ( E) \right).$$

Recall that we have already chosen $\varepsilon_0$ previously.
It remain to choose $\delta>0$. For a given $\varepsilon$, let
  $\delta$ be small enough to ensure that when $\omega^{_{V_d}} ( E) / \omega^{_{V_d}} (
  \Delta_r) < \delta$, then the length of good $\varepsilon_0$-cover is sufficiently large so that $(S^r)^2(u) \gtrsim K $ and thus $\sigma ( E) / \sigma ( \Delta)
  \lesssim K^{-1} < \varepsilon$. This concludes our proof.\qed


\section{Proof of Theorem \ref{T:Main4}}

This proof of Theorem \ref{T:Main4} is based on the following lemma from \cite{DH}.

\begin{lemma}\label{L:Square} (Lemma 3.3 of \cite{DH}) Let $\Omega$ be an admissible domain from Definition \ref{D:domain} of character $(\ell,N,C_0)$. Let $L=\partial_t-\di(A\nabla\cdot)$ be a parabolic operator with matrix $A$ satisfying uniform ellipticity with constants $\lambda$ and $\Lambda$, and such that, for all $0<r\le r',$

\begin{equation}\label{carlIII}
d\mu=\delta(X)\left(\sup_{B_{\delta(X)/2}(X)}|\nabla a_{ij}|\right)^2\,dX\,dt
\end{equation}
is a Carleson measure with norm $\|\mu\|_{Carl}$.

Then there exist a constant $C=C(\lambda,\Lambda,N,C_0)$  such that
for any solution $u$ with boundary data $f$ on any ball $\Delta_r\subset\partial\Omega$ with $r\le \min\{r'/4, r_0/(4C_0)\}$ (c.f. Definition \ref{D:domain} for the meaning of $r_0$ and $C_0$) we have

\begin{equation}\label{E:sq10aa}
\int_{T(\Delta_r)} |\nabla u |^2 x_0\,dX\,dt
 \leq     C(1+\|\mu\|_{C,2r_0})(1+\ell^2) \int_{\Delta_{2r}} (N^{2r})^{2}(u) \,dX\,dt.
\end{equation}
Here $N^{2r}$ denotes the truncated non-tangential maximal function.
\end{lemma}

\noindent{\it Remark.} Let $u$ be a solution of $Lu=0$ in $\Omega$ with bounded boundary data $f$. Since by the maximum principle $\|u\|_{L^\infty}\le\|f\|_{L^\infty}$ and $\|N^{2r}\|^2_{L^\infty(\Delta_{2r})}\le \|u\|^2_{L^\infty(\Omega)}\sigma(\Delta_{2r})$,
it follows from \eqref{E:sq10aa} that
$$\int_{T(\Delta_r)} |\nabla u |^2 x_0\,dX\,dt
 \leq     C(1+\|\mu\|_{C,2r_0})(1+\ell^2) \|f\|_{L^\infty}^2\sigma(\Delta_{2r}),$$
hence by doubling for all $\Delta_r$ with $0<r\le r''$ we have
\begin{equation}\label{E:sq10ab}
\sigma^{-1}(\Delta_r)\int_{T(\Delta_r)} |\nabla u |^2 x_0\,dX\,dt
 \leq     C(1+\|\mu\|_{C,2r_0})(1+\ell^2) \|f\|^2_{L^\infty},
\end{equation}
which by Theorem \ref{T:Main2} shows $A_\infty$. We shall track how the $L^p$ solvability depends on the ellipticity $\lambda,\Lambda$
and the constant $K=C(1+\|\mu\|_{C,2r_0})(1+\ell^2)$ in the estimate \eqref{E:sq10ab}.

We have the following:

\begin{lemma}\label{L:trackP} Let $\Omega$ be an admissible domain from Definition \ref{D:domain} and $L=\partial_t-\di(A\nabla\cdot)$ be a parabolic operator with matrix $A$ satisfying uniform ellipticity with constants $\lambda$ and $\Lambda$. Suppose that for all solutions $u$ with boundary data $f \in L^{\infty}$ we have
\begin{equation}\label{E:sq10ac}
\sigma^{-1}(\Delta_r)\int_{T(\Delta_r)} |\nabla u |^2 x_0\,dX\,dt
 \leq     K \|f\|^2_{L^\infty},
\end{equation}
for all $0<r\le r''$. Then there exists $p_0=p_0(\lambda,\Lambda,K)>0$ such that for all $p_0<p \leq \infty$ the $L^p$ Dirichlet problem is solvable for the operator $L$.
\end{lemma}

\noindent{\it Proof.} We revisit the proof of Theorem \ref{T:Main2} from the previous section, tracking how the the result depends on various parameters.
Assuming \eqref{E:sq10ac} we have established that for $E\subset\Delta_r\subset \Delta_d$
\begin{equation}
(\beta,\varepsilon)_{A_\infty}:\quad\frac{\omega^{V_d}(E)}{\omega^{V_d}(\Delta_r)}<\beta\quad\Longrightarrow\quad \frac{\sigma(E)}{\sigma(\Delta_r)}<\varepsilon\quad\mbox{for}\quad \varepsilon=\frac{C(\lambda,\Lambda,\ell)(1+K)}{-\log \beta}.\label{E:ainfty}
\end{equation}
Here we took into account that the $\varepsilon_0$ in the good-$\varepsilon_0$ cover depended on $\lambda,\,\Lambda$ and $\ell$ and that the length of
the the good-$\varepsilon_0$ cover was $\approx -\varepsilon_0\log\beta$. We have established this for balls on $\partial\Omega$ but same will also hold
for parabolic cubes (in fact the metric $d$ can be defined in such that \lq\lq balls" in $d$ are just parabolic cubes).

Given the $(\beta,\varepsilon)_{A_\infty}$ statement we want to show that $(\beta,\varepsilon)_{A_\infty}\Longrightarrow (\beta,\alpha)_{A'_\infty}$,
for $\alpha=1-\varepsilon$ where
$$(\beta,\alpha)_{A'_\infty}:\qquad\sigma\left(\{(x,t)\in \Delta:\,K^{V_d}(x,t)>\beta\omega^{V_d}(\Delta)/\sigma(\Delta)\}\right)>\alpha\sigma(\Delta),\quad\mbox{for all }\Delta\subset\Delta_d.$$
Here $K^{V_d}$ is the Radon-Nykodim derivative $\frac{d\omega^{V_d}}{d\sigma}$. To see this let
$$E=\{(x,t)\in\Delta:\,K^{V_d}(x,t)\le\beta\omega^{V_d}(\Delta)/\sigma(\Delta)\}.$$
Then
$$\omega^{V_d}(E)=\int_E K^{V_d}\,d\sigma \le \frac{\beta\omega^{V_d}(\Delta)}{\sigma(\Delta)}\int_E\, d\sigma\le \beta\omega^{V_d}(\Delta).$$
Hence since $\frac{\omega^{V_d}(E)}{\omega^{V_d}(\Delta)}<\beta$ by $(\beta,\varepsilon)_{A_\infty}$ we have that $\frac{\sigma(E)}{\sigma(\Delta)}<\varepsilon$. Hence
$$\sigma(E^c)=\sigma(\Delta)-\sigma(E)>(1-\varepsilon)\sigma(\Delta),$$
which is exactly $(\beta,1-\varepsilon)_{A'_\infty}$. We use standard arguments to show that $(\beta,\alpha)_{A'_\infty}$ implies a reverse-H\"older inequality for $K^{V_d}$.
That is, $(\beta,\alpha)_{A'_\infty}$ implies that $\omega^{V_d}$ belongs to $B_{1+\delta}(\sigma)$, where we track the dependence of $\delta$ on $n, \varepsilon, \beta.$


As in \cite{CF}, for any $\lambda>m_{\Delta}:=\omega^{V_d}(\Delta)/\sigma(\Delta)$ the Calder\'on-Zygmund lemma produces a
family of pairwise disjoint cubes $Q_i$ such that $K^{V_d}\le \lambda$ for a.e. $x\in \Delta\setminus \bigcup_i Q_i$ and
$$\lambda\le \sigma(Q_i)^{-1}\int_{Q_i}K^{V_d}d\sigma\le 2^n \lambda.$$
Hence by $(\beta,\alpha)_{A'_\infty}$ we obtain
\begin{eqnarray}
&&\int_{\{(x,t)\in\Delta:\,K^{V_d}(x,t)>\lambda\}}K^{V_d}d\sigma\le \sum_i \int_{Q_i}K^{V_d}d\sigma \le 2^n\lambda \sum_i \sigma(Q_i)\nonumber\\
&\le& \frac{2^n\lambda}{\alpha}\sum_i \sigma\left(\{(x,t)\in Q_i:\,K^{V_d}(x,t)>\beta\omega^{V_d}(Q_i)/\sigma(Q_i)\}\right)\nonumber\\
&\le& \frac{2^n\lambda}{\alpha}\sum_i \sigma\left(\{(x,t)\in Q_i:\,K^{V_d}(x,t)>\beta\lambda\}\right)\nonumber\\
&\le& \frac{2^n\lambda}{\alpha}\sigma\left(\{(x,t)\in \Delta:\,K^{V_d}(x,t)>\beta\lambda\}\right)\nonumber
\end{eqnarray}

Hence using this we get for the integral:
$$\int_{m_\Delta}^\infty \lambda^{\delta-1}\int_{\{(x,t)\in\Delta:\,K^{V_d}(x,t)>\lambda\}}K^{V_d}d\sigma d\lambda\le \frac{2^n}{\alpha}\int_0^\infty
\lambda^\delta\sigma\left(\{(x,t)\in \Delta:\,K^{V_d}(x,t)>\beta\lambda\}\right)d\lambda$$
which further (after a substitution $t=\beta\lambda$) equals to
$$=\frac{2^n}{\alpha\beta^{1+\delta}}\int_0^\infty t^\delta \sigma\left(\{(x,t)\in \Delta:\,K^{V_d}(x,t)>t\}\right)dt=\frac{2^n}{\alpha\beta^{1+\delta}(1+\delta)}\int_{\Delta}(K^{V_d})^{1+\delta}d\sigma.$$
On the other hand by Fubini
\begin{eqnarray}
&&\int_{m_\Delta}^\infty \lambda^{\delta-1}\int_{\{(x,t)\in\Delta:\,K^{V_d}(x,t)>\lambda\}}K^{V_d}d\sigma d\lambda\nonumber\\
&\ge& \int_{\{(x,t)\in\Delta:\,K^{V_d}(x,t)>m_\Delta\}}K^{V_d}(x,t)\int_{m_{\Delta}}^{K^{V_d}(x,t)}\lambda^{\delta-1}d\lambda\,d\sigma\nonumber\\
&=&\int_{\{(x,t)\in\Delta:\,K^{V_d}(x,t)>m_\Delta\}}K^{V_d}(x,t)\left[\frac{K^{V_d}(x,t)^\delta}{\delta}-\frac{m_\Delta^\delta}{\delta}\right]d\sigma\nonumber\\
&\ge&\frac1\delta\int_{\Delta}(K^{V_d})^{1+\delta}d\sigma-\frac{m_\Delta^{1+\delta}}{\delta}\sigma(\Delta).\nonumber
\end{eqnarray}
It follows that
$$\left(\frac1\delta-\frac{2^n}{\alpha \beta^{1+\delta}(1+\delta)}\right)\sigma(\Delta)^{-1}\int_{\Delta}(K^{V_d})^{1+\delta}d\sigma\le \frac1\delta \left(\sigma(\Delta)^{-1}\int_{\Delta}K^{V_d}d\sigma\right)^{1+\delta},$$
from which our claim follows.\vglue1mm


Thus if $(\beta,\varepsilon)_{A_\infty}$ holds, then $\omega^{V_d} \in B_{1+\delta}(\sigma)$ for all $0<\delta<\delta_0$ where
$$\frac1{\delta_0}=\frac{2^n}{(1-\varepsilon) \beta^{1+\delta_0}(1+\delta_0)}.$$
The duality relationship tells us that the $L^p$ Dirichlet problem is solvable for $p=(1+\delta)/\delta$. To obtain an estimate from below on $p$ we may assume that $\delta<1$.
Hence if 
$$p=\frac{1+\delta}{\delta}> \frac{2^{n}}{(1-\varepsilon)\beta^2}\ge \frac{2^n}{(1-\varepsilon) \beta^{1+\delta}},$$
and if we choose $\varepsilon=1/2$, and the corresponding $\beta=\exp(-2C(\lambda,\Lambda,\ell)(1+K))$ given by \eqref{E:ainfty} then for
\begin{equation}
p>p_0:=2^{n+1}\exp(4C(\lambda,\Lambda,\ell)(1+K))\label{Estim_p}
\end{equation}
the $L^p$ Dirichlet problem for the operator $L$ satisfying \eqref{E:sq10ac} is solvable. This establishes Lemma \ref{L:trackP}.\qed

In particular by \eqref{E:sq10ab}, Lemma \ref{L:trackP} applies directly to operators $L=\partial_t-\di(A\nabla\cdot)$ on admissible domains $\Omega$ whose matrix $A$ satisfies
$$d\mu=\delta(X)\left(\sup_{B_{\delta(X)/2}(X)}|\nabla a_{ij}|\right)^2\,dX\,dt$$
is a Carleson measure. Thus, there exists $p_0=p_0(n,\lambda,\Lambda,\|\mu\|_{Carl},L)<\infty$ such that the $L^p$ Dirichlet problem for such an operator $L$
is solvable for all $p>p_0$. In particular this implies that Theorem \ref{T:Main4} holds for such operators $L$. Indeed, by \cite{DH} for all $p\ge 2$ if $\max\{\ell^2, \|\mu\|_{Carl}\}$ is sufficiently small then the $L^p$
Dirichlet problem is solvable for an operator $L$. Let $C(p)=C(p,\lambda,\Lambda,n)>0$ be the largest number for which  the condition $\max\{\ell^2, \|\mu\|_{Carl}\}<C(p)$ implies $L^p$ solvability. To show that $C(p)\to\infty$ we only have to prove two claims. The first one is that the function $C(p)$ is monotone non-decreasing in $p$. This is due to the fact that $L^p$ solvability implies $L^{q}$ solvability for all $q>p$.
The second that is that for an arbitrary fixed $M>0$, if $\max\{\ell^2, \|\mu\|_{Carl}\}<M$ then there exists
$p<\infty$ such that $C(p)\ge M$. An estimate of how large such a $p$ must be is given by \eqref{Estim_p}. Observe that $\ell$ and $K$ on the righthand side of \eqref{Estim_p} are both bounded by $M$, and hence the value of $p$ for which $C(p)\ge M$ only depends on $n,\,\lambda,\,\Lambda$ and $M$. Combining these two claim yields
$$\lim_{p\to\infty}C(p)=\infty,$$
as desired.\medskip

We will use the theory of perturbation of operators to conclude that Theorem \ref{T:Main4} holds for operators, where the condition on the gradient has been replaced by the oscillation condition \eqref{carlII}.  
Let $L_0 = \partial_t -\di (A_0 \nabla \cdot)$ be an operator satisfying \eqref{carlII} with Carleson norm $K_0$.

We will proceed in two steps, introducing two intermediate operators $L_1$ and $L_2$ to which $L_0$ will be compared.

Following \cite{DPP}, create a new operator $L_2$, namely 
$L_2  = \partial_t-\di (A_2 \nabla \cdot)$, where $A_2$ is the mollification of $A_0$ obtained by convolving the coefficients with a smooth bump function. Then  the coefficients of $L_2$ satisfy the Carleson gradient condition \eqref{carlIII} with norm $K_2:=\|d\widetilde{\mu}\|_{Carl}$, bounded by a multiple of $K_0$.
(See the proof of  Corollary 2.3 of \cite{DPP} and the proof of Theorem 3.1 in
\cite{DH} for more details on this construction.)




Precisely, the difference between the coefficients of $L_0$ and $L_2$ satisfies the perturbation-Carleson condition (\cite{N}, \cite{Sw}, \cite{FKP}) with constant $C(K_0)$, a multiple of $K_0$:
\begin{equation}\label{carlperturb}
  \sup_{\Delta_r}  \,\,  (\sigma(\Delta_r))^{-1} \int_{\Delta_r}  \int_{\Gamma^r} (\delta(X)^{-2-n}\left(\sup_{B_{\delta(X)/2}(X)}|A_0(X) - {A}_2(X)|\right)^2\, dXdt\, d\sigma <  C(K_0)
  \end{equation}
  
Here, $\Gamma^r$ is the truncated cone of \eqref{D:Gammagr} and, after integrating, the condition \eqref{carlperturb} is equivalent to stating that 
$\delta(X)^{-2-n}\left(\sup_{B_{\delta(X)/2}(X)}|A_0(X) - A_2(X)|\right)^2\, dXdt$ is itself a Carleson measure. 

However, it is be useful to write this in the form written in \eqref{D:Gammagr}, to recall that the Carleson measure condition tells us that the (truncated) area 
integral is bounded on a large fraction of $\Delta_r$. 
We now, as in  \cite{N}, \cite{FKP}, introduce another operator $L_1 = \partial_t -\di (A_1 \nabla \cdot)$, which will have the stronger property that 

\begin{equation}\label{finitesquare}  A_r^2(Q): = \int_{\Gamma^r(Q)} (\delta(X)^{-2-n}\left(\sup_{B_{\delta(X)/2}(X)}|A_0(X) - A_1(X)|\right)^2\, dXdt  \leq C_1(K_0)
\end{equation}

The construction of $A_1$ proceeds as follows.  Fix a set $F \subset \Delta_r$ with $\sigma(F \cap \Delta_r) >  \sigma(\Delta_r)/2$ and $A_r^2(Q) < CK_0$ on $F$.
As in \cite{N}, Section 3, a sawtooth region is formed over $F$: the new matrix $A_1$ will equal $A_2$ in that sawtooth region over $F$, and equal $A_0$ otherwise. It is argued in \cite{N}, following \cite{FKP}, that
the resulting operator $L_1$ satisfies \eqref{finitesquare}.

\medskip

We summarize the steps, and make some comments regarding tracking the dependence on $p$ for solvability of the $L^p$ Dirichlet problem.

\smallskip

{\it Step 1.}  The solvability of $D_{p}$ for $L_2$, for all $p > p_2$, with an estimate of the dependence of $p_2$ on $K_2$ was carried out above in the proof of Lemma \ref{L:trackP}. Thus
(suppressing the dependence on the corkscrew point), we have that the $(\beta, \epsilon)_{A_\infty}$ condition holds for $\omega_{L_2}$, which gives a 
constant $\delta = \delta(\beta, \epsilon)$ such that $\omega_{L_2} $ belongs to $B_{1+\delta}$.  

\smallskip

{\it Step 2.}  The solvability of $D_{p}$ for $L_1$, for all $p > p_1$, and with $p_1 = p_1(\beta,\epsilon)$ is a consequence of the construction of sawtooth regions and a comparison of the parabolic measures
$d\omega_{L_2}$ and $d\omega_{L_1}$.  The constants will depend on domain parameters, and will introduce no further dependence upon the Carleson norm. The comparison of these two measures is carried out in \cite{N}, following the construction in \cite{DJK} in the elliptic case.

\smallskip

{\it Step 3.}  The solvability of a $D_{p_0}$ for $L_0$, will result from the a chain of comparisons starting with $L_1$ and ending with $L_0$. The parameter $p_0$ can be tracked explicitly through the transitivity of the reverse H\"older classes, $B_q$.
The method (\cite{FKP}) is as follows: Form the family of parabolic operators $L_{s}$, moving from $L_1$ at $s=1$ to $L_0$ at $s=0$ where 
$A_s = (1-s)A_0 + sA_1$  and $L_s = \partial_t-\di (A_s \nabla \cdot)$.
Theorem 6.6 of \cite{N} provides a small $\epsilon_0$ that depends only domain parameters and ellipticity for which $\omega_{L_{k\eta}}$ belongs to $B_2$ with respect to $\omega_{L_{(k+1)\eta}}$, for any $\eta < \epsilon_0$.

\noindent {\it Remark.} Suppose that $\omega_0$, $\omega_1$ and $\omega_2$ are weights satisfying $\omega_1\in B_{p_1}(\omega_0)$ and $\omega_2 \in B_{p_2}(\omega_1)$ with constants
$\|\omega_1\|_{B_{p_1}(\omega_0)}$ and $\|\omega_2\|_{B_{p_2}(\omega_1)}$ respectively.
Then
$\omega_2\in B_r(\omega_0)$ where $r = \frac{p_1p_2}{p_1+p_2-1}$ and the $\|\omega_2\|_{B_r(\omega_0)}$ depends upon
$p_1,p_2, \|\omega_1\|_{B_ p(\omega_0)}, \|\omega_2\|_{B_ p(\omega_1)}$.\medskip

\noindent  To conclude, this will be applied approximately $\epsilon_0^{-1}$ times and we find that $\omega_{L_0}\in B_{1+\delta'}(d\sigma)$ for some positive $\delta'=\delta'(n,\lambda,\Lambda,\ell,\|\mu\|_{Carl})$ from which $L^p$ solvability of the Dirichlet problem follows for all $p>(1+\delta')/\delta'$. As above, this implies $C(p)$ in Theorem \ref{T:Main4} has the property that $C(p)\to\infty$ as $p\to\infty$.

\end{document}